\documentclass{article}
\usepackage{amsmath,amssymb,enumerate,epsfig,bbm,calc,theorem,ifthen,capt-of}

\newcommand{\mathd}{\mathrm{d}}
\newcommand{\tmmathbf}[1]{\ensuremath{\boldsymbol{#1}}}
\newcommand{\tmop}[1]{\ensuremath{\operatorname{#1}}}
\newcommand{\tmstrong}[1]{\textbf{#1}}
\newcommand{\tmtextit}[1]{{\itshape{#1}}}
\newcommand{\tmtextup}[1]{{\upshape{#1}}}
\newenvironment{enumeratenumeric}{\begin{enumerate}[1.] }{\end{enumerate}}

\newenvironment{proof}{\noindent\textbf{Proof\ }}{\hspace*{\fill}$\Box$\medskip}
\newtheorem{lemma}{Lemma}
{\theorembodyfont{\rmfamily}\newtheorem{remark}{Remark}}
\newtheorem{theorem}{Theorem}
\newcommand{\tmfloatcontents}{}
\newlength{\tmfloatwidth}
\newcommand{\tmfloat}[5]{
  \renewcommand{\tmfloatcontents}{#4}
  \setlength{\tmfloatwidth}{\widthof{\tmfloatcontents}+1in}
  \ifthenelse{\equal{#2}{small}}
    {\ifthenelse{\lengthtest{\tmfloatwidth > \linewidth}}
      {\setlength{\tmfloatwidth}{\linewidth}}{}}
    {\setlength{\tmfloatwidth}{\linewidth}}
  \begin{minipage}[#1]{\tmfloatwidth}
    \begin{center}
      \tmfloatcontents
      \captionof{#3}{#5}
    \end{center}
  \end{minipage}}

\begin{document}

\title{Level Set Dynamics and the Non-blowup of the 2D Quasi-geostrophic
Equation} 

\author{J. Deng\thanks{School of Math. Sciences, Fudan University,
Shanghai, China. Email: jdeng@fudan.edu.cn}, 
T. Y. Hou\thanks{Applied and Comput. Math, Caltech, Pasadena,
CA 91125. Email: hou@acm.caltech.edu.}, 
R. Li\thanks{Applied and Comput. Math, Caltech, Pasadena,
CA 91125. Email: rli@acm.caltech.edu.}, 
X. Yu\thanks{Dept. of Math., UCLA,
Los Angeles, CA 90095. Email: xinweiyu@math.ucla.edu.}} 

\maketitle

\begin{abstract}
  In this article we apply the technique proposed in Deng-Hou-Yu 
  {\cite{deng-hou-yu:2005}} to study the level set dynamics of the 
  2D quasi-geostrophic equation. Under certain assumptions on the
  local geometric regularity of the level sets of $\theta$, we obtain
  global regularity results with improved growth estimate on 
  $\left| \nabla^{\bot} \theta \right|$. We further perform numerical 
  simulations to study the local geometric properties of the level sets
  near the region of maximum $\left| \nabla^{\bot} \theta \right|$.
  The numerical results indicate that the assumptions on the local
  geometric regularity of the level sets of $\theta$ in our theorems 
  are satisfied. Therefore these theorems provide a good explanation
  of the double exponential growth of $\left| \nabla^{\bot} \theta \right|$ 
  observed in this and past numerical simulations.
\end{abstract}

\section{\label{sec:intro}Introduction}

The study of global existence/finite-time blow-up of the two-dimensional
quasi-geostrophic (subsequently referred to as 2D QG for simplicity ) equation
has been an active research area in the past ten years, partly due to its close
connection to the 3D incompressible Euler equations (Constantin-Majda-Tabak
{\cite{constantin-majda-tabak:1994}}, Cordoba {\cite{cordoba:1998}},
Cordoba-Fefferman {\cite{cordoba-fefferman:2002}}). The 2D QG equation has its
origin in modeling rotating fluids on the earth surface (Pedlosky
{\cite{pedlosky:1987}}). The equation describes the transportation of a scalar
quantity $\theta$:
\begin{equation}
  \label{eq:qg} D_t \theta \equiv \theta_t + u \cdot \nabla \theta = 0
\end{equation}
with initial conditions $\theta \mid_{t = 0} = \theta_0$. The relation between
$\theta$ and the velocity $u$ is given by
\begin{equation}
  \label{eq:qg.u.psi} u = \nabla^{\bot} \psi, \hspace{1em} \psi = \left( -\bigtriangleup
  \right)^{- \frac{1}{2}} \left( - \theta \right)
\end{equation}
where
\begin{equation}
  \label{eq:qg.psi} \nabla^{\bot} \psi \equiv \left( - \frac{\partial
  \psi}{\partial x_2}, \frac{\partial \psi}{\partial x_1} \right)^T
\end{equation}
and
\begin{equation}
  \left( -\bigtriangleup \right)^{- \frac{1}{2}} \psi \equiv \int e^{2 \pi ix \cdot k}
  \frac{1}{2 \pi \left| k \right|} \hat{\psi} \left( k \right) \mathd k
\end{equation}
where $\hat{\psi} \left( k \right) = \int e^{-2 \pi ix \cdot k} \psi 
\left( x \right) \mathd x$ is the Fourier transform of $\psi \left( x \right)$.

As pointed out by Constantin-Majda-Tabak
{\cite{constantin-majda-tabak:1994}}, the 2D QG equation bears striking
mathematical and physical analogy to the 3D incompressible Euler equations.
They both exhibit similar geometric/analytic structures. 
In particular, one can derive a necessary and sufficient blow-up condition 
for the 2D QG equation similar to the well-known Beale-Kato-Majda
criterion (Beale-Kato-Majda {\cite{beale-kato-majda:1984}}).
More precisely, the solution to the 2D QG equation (\ref{eq:qg})
becomes singular at time $T_{\ast}$ if and only if
\begin{equation}
  \label{eq:qg.BKM} \int_0^{T_{\ast}} \left\| \nabla^{\bot} \theta \left(
  \cdot, t \right) \right\|_{L^{\infty}} \mathd t = + \infty .
\end{equation}
Thus, $ \nabla^{\bot} \theta $ plays a role similar to the vorticity $\omega$ in
the 3D Euler equations. 
Furthermore, as in the 3D incompressible Euler equations, the velocity
$u$ is related to $\nabla^{\bot} \theta$ by an order $-1$ singular 
integral operator.

On the other hand, in some aspects the 2D QG equation behaves much better
than the 3D incompressible Euler equations. For example, it was shown in
Cordoba {\cite{cordoba:1998}} that $\left\| u \left( \cdot, t \right)
\right\|_{L^{\infty}}$ is bounded by $\log \left\| \nabla^{\bot} \theta \left(
\cdot, t \right) \right\|_{L^{\infty}}$ for any time $t$, while for the 3D
incompressible Euler equations, $\left\| u \left( \cdot, t \right)
\right\|_{L^{\infty}}$ may grow as fast as $\left\| \omega \left( \cdot, t
\right) \right\|_{L^{\infty}}^{1 / 2}$ according to Kelvin's circulation
theorem.

Much effort has been made to obtain global existence for the 2D QG equation.
In Constantin-Majda-Tabak {\cite{constantin-majda-tabak:1994}}, it was shown
that if the direction field $\xi \equiv \nabla^{\bot} \theta / \left|
\nabla^{\bot} \theta \right|$ remains smooth in a region, then no finite-time
singularity is possible in that region. In particular, if this region contains
maximum $\left| \nabla^{\bot} \theta \right|$ all the time, then the solution
remains regular global in time. Based on this understanding, they conjectured
that the 2D QG equation with initial level sets of the hyperbolic saddle type
is likely to develop a finite time singularity. They further presented 
some numerical evidence which supports a finite time singularity for the
2D QG equation. Later on, Ohkitani-Yamada {\cite{ohkitani-yamada:1997}},
Constantin-Nie-Schorghofer {\cite{constantin-nie-schorghofer:1998}},
{\cite{constantin-nie-schorghofer:1999}} re-did the numerical simulations with
higher resolutions, and revealed that the growth of $\left\| \nabla^{\bot}
\theta \right\|_{L^{\infty}}$ in time is no faster than double exponential. 
Around the same time, Cordoba {\cite{cordoba:1998}} proved that
the growth of $\left\| \nabla^{\bot} \theta \right\|_{L^{\infty}}$ in 
the hyperbolic saddle scenario is bounded by quadruple exponential 
under the assumption that near the saddle point the level sets of 
$\theta$ exhibits certain self-similar structure. This result was 
subsequently simpified by Cordoba-Fefferman in {\cite{cordoba-fefferman:2002}}.

In this paper, we take a different approach first proposed by Deng-Hou-Yu
{\cite{deng-hou-yu:2005}} to study the singularity problem of the 3D Euler
equations. Following the similar approach for the 2D QG equation,
we study the Lagrangian evolution of some
localized segments of level sets carrying large $\left| \nabla^{\bot} \theta
\right|$. By exploring the incompressibility condition of $\nabla^{\bot}
\theta$ and using the local geometric properties of level sets, we obtain
estimates for the growth of $\left| \nabla^{\bot} \theta \right|$ by 
studying the stretching of such level set segments. We find that, when there
is one level set segment of length $O \left( \frac{1}{\tmop{loglog} \left\|
\nabla \theta \right\|_{L^{\infty}}} \right)$ on which the maximum $\left|
\nabla^{\bot} \theta \right|$ is comparable to the global maximum, and along
which $\nabla \left( \frac{\nabla^{\bot} \theta}{\left| \nabla^{\bot} \theta
\right|} \right)$ is bounded by $O \left( \tmop{loglog} \left\| \nabla \theta
\right\|_{L^{\infty}} \right)$, the growth rate of $\left\| \nabla^{\bot}
\theta \right\|_{L^{\infty}}$ is bounded by triple exponential. In particular,
when the length of the segment and the bound of $\nabla \left(
\frac{\nabla^{\bot} \theta}{\left| \nabla^{\bot} \theta \right|} \right)$ are
both $O \left( 1 \right)$, we can improve our estimate on the growth
of $\left\| \nabla^{\bot} \theta \right\|_{L^{\infty}}$ and bound it
by double exponential. The double exponential estimate is sharp according to 
recent numerical simulations (Ohkitani-Yamada
{\cite{ohkitani-yamada:1997}}, Constantin-Nie-Schorghofer
{\cite{constantin-nie-schorghofer:1998}},
{\cite{constantin-nie-schorghofer:1999}}). 

We also perform
careful numerical experiments to study the local geometric properties 
of the level sets in a region containing maximum $\left|\nabla^{\bot} \theta\right|$. 
Our numerical results indicate that this region of large
$\left|\nabla^{\bot} \theta\right|$ and the region of large
$\left|\nabla \left( \frac{\nabla^{\bot} \theta}{\left| \nabla^{\bot} \theta
\right|} \right)\right|$ are essentially disjoint. Furthermore, there exists
an $O \left( 1 \right)$ level set segment within this region of 
large $\left|\nabla^{\bot} \theta\right|$ along which 
$\nabla \left( \frac{\nabla^{\bot} \theta}{\left| \nabla^{\bot} \theta
\right|} \right)$ is bounded. Thus our second theorem applies, which
implies that $\left\|\nabla^{\bot} \theta(\cdot,t)\right\|_{L^\infty}$ is bounded by double exponential
in time. In some sense, our theoretical results capture the essential
feature of the dynamic growth of the 2D QG equation.

The rest of the paper is organized as follows. In Section \ref{sec:main-results} 
we give an overview of the main results. In Section \ref{sec:numerical} we present
numerical results which illustrate the local geometric properties of the
level sets in the region containing maximum $\left|\nabla^{\bot} \theta\right|$.
Finally, in
Section \ref{sec:key-est} we estimate level set stretching and obtain an
estimate for the growth of $\left| \nabla^{\bot} \theta \right|$, and 
prove the main theorems.

\section{\label{sec:main-results}Main Results}

We present the main results in this section. Denoting $\left\| \nabla^{\bot}
\theta \left( \cdot, t \right) \right\|_{L^{\infty}}$ by $\Omega \left( t
\right)$, we consider, at time $t$, a level set segment $L_t$ along which the
maximum of $\left| \nabla^{\bot} \theta \right|$ (denoted by $\Omega_L \left( t
\right)$ in the following) is comparable to $\Omega \left( t \right)$. Denote
by $L \left( t \right)$ the arc length of $L_t$, $\xi$ the tangential, and
$\tmmathbf{n}$ the normal unit vector of $L_t$. The direction of $\xi$ and
$\tmmathbf{n}$ are determined as follows: $\xi = \frac{\nabla^{\bot}
\theta}{\left| \nabla^{\bot} \theta \right|}$, $\tmmathbf{n}= \frac{\xi \cdot
\nabla \xi}{\left| \xi \cdot \nabla \xi \right|}$. We define $U_{\xi} \left( t
\right) = \max_{x, y \in L_t} \left| \left( u \cdot \xi \right) \left( x, t
\right) - \left( u \cdot \xi \right) \left( y, t \right) \right|$, $U_n \left(
t \right) = \max_{x \in L_t} \left| \left( u \cdot \tmmathbf{n} \right) \left(
x, t \right) \right|$, $M \left( t \right) = \max_{x \in L_t} \left| \nabla
\cdot \xi \right|$, and $K \left( t \right) = \max_{x \in L_t} \kappa$ where
$\kappa = \left| \xi \cdot \nabla \xi \right|$ is the (unsigned) curvature. We
should point out that our theorems only requires $L_t$ to be a subset of $X
\left( L_{t'}, t', t \right)$, the flow image of $L_{t'}$ at time $t$, for $t'
< t$. With these notations, we present our main results.

\begin{theorem}
  \label{thm:triple}Assume that there is a family of level set segments $L_t$
  and $T_0 \in \left[ 0, T_{\ast} \right)$ such that $X \left( L_{t_0}, t_0, t
  \right) \supseteq L_t$ for all $T_0 \leqslant t_0 < t < T_{\ast}$. Also
  assume that $\Omega \left( t \right)$ is monotonically increasing and
  $\Omega_L \left( t \right) \geqslant c_0 \Omega \left( t \right)$ for some
  $0 < c_0 \leqslant 1$ for all $t \in \left[ T_0, T_{\ast} \right)$. Then the
  classical solution of the 2D QG equation can be extended beyond $T_{\ast}$
  as long as the following conditions are satisfied: there exists constants
  $c_L, C_0 > 0$, such that
  \begin{itemize}
    \item (H1). $L \left( t \right) \geqslant \frac{c_L}{\tmop{loglog} \Omega
    \left( t \right)}$, and
    
    \item (H2). $M \left( t \right) L \left( t \right), \hspace{1em} K \left(
    t \right) L \left( t \right) \leqslant C_0$.
  \end{itemize}
  Furthermore, for $t \in \left[ 0, T_{\ast} \right)$, we have the following
  triple exponential estimate:
  \begin{equation}
    \label{eq:triple.est} \tmop{logloglog} \Omega \left( t \right) \leqslant
    C_1 t + C_2
  \end{equation}
  for some constants $C_1, C_2 > 0$ independent of $t$. 
\end{theorem}

\begin{remark}
  If we further assume that 
  $\Omega (t) \leqslant (T_{\ast} - t)^{- B}$ for some $B<+\infty$, 
  we can easily prove non-blowup of the solution with condition (H1) replaced by 
  $L \left( t \right) \geqslant c_L \left( T_{\ast} - t \right)^A$ 
  for any $A < 1$. See Yu {\cite{yu:2005}} for details. 
\end{remark}

With stronger assumptions on the regularity of 
$M \left( t \right), L \left( t \right), K \left( t \right)$, we can 
obtain a sharper growth estimate, which yields double exponential growth 
of $\Omega \left( t \right)$. This growth rate is consistent with the
observations in recent numerical simulations 
(Ohkitani-Yamada {\cite{ohkitani-yamada:1997}}, 
Constantin-Nie-Schorghofer {\cite{constantin-nie-schorghofer:1998}},
{\cite{constantin-nie-schorghofer:1999}}).

\begin{theorem}
  \label{thm:double}Assume that all the assumptions and conditions in Theorem
  \ref{thm:triple} hold, except that (H1) is replaced by
  \begin{itemize}
    \item (H1'). $L \left( t \right) \geqslant c_L$.
  \end{itemize}
  Then the estimate (\ref{eq:triple.est}) can be improved to
  \[ \tmop{loglog} \Omega \left( t \right) \leqslant C_1' t + C_2' \]
  for some constants $C_1'$ and $C_2'$ independent of time.
\end{theorem}

\begin{remark}
  In Section \ref{sec:numerical}, we will perform numerical experiments
  to study the local geometric properties of the level sets. In particular,
  we will show that the conditions (H1') and (H2) are consistent
  with our numerical results. Furthermore, for certain hyperbolic saddle 
  scenario similar to the one studied by Cordoba in {\cite{cordoba:1998}}
  but with additional assumptions, we can prove that the conditions (H1') 
  and (H2) are actually satisfied. This gives a partial theoretical justification
  of the conditions (H1') and (H2) in certain hyperbolic saddle scenario. 
  The key observation is that the maximum 
  $\left| \nabla^{\bot} \theta \right|$ along any level set is located 
  away from the saddle point. Therefore although $\left|\nabla \xi\right|$ is large 
  near the saddle point, it is bounded along level set segments considered 
  in Theorem \ref{thm:double}. For details, see Yu {\cite{yu:2005}}.
\end{remark}

\section{\label{sec:numerical}Numerical Results}

To further understand the dynamics of the 2D QG system
(\ref{eq:qg})--(\ref{eq:qg.psi}) and to study the local geometric
properties of the level sets, we perform careful numerical
simulations of the 2D QG equation. Specifically, we would like
to track the dynamic evolution of the following two regions
which characterize the geometric regularity of the level sets
around the points of maximum $\left| \nabla^{\bot} \theta \right|$:
\begin{enumeratenumeric}
  \item The region of large $\left| \nabla^{\bot} \theta \right|$, and
  
  \item The region of large $\left| \nabla \xi \right|$.
\end{enumeratenumeric}
As we will demonstrate from our numerical simulations, these two regions 
are essentially disjoint. Although they both undergo severe stretching
and thining as the flow evolves, their intersection occupies just a small 
portion of either one of these two regions. 
This numerical result is quite surprising. 
It shows that in the very localized region where the maximum of
$\left| \nabla^{\bot} \theta \right|$ is attained, the level sets across
this region are regular. On the other hand, in the region where
$\left| \nabla \xi \right|$ is very large, the value of 
$\left| \nabla^{\bot} \theta \right|$ is relatively small compared with
its global maximum. This complementary local geometric regularity of
level set filaments seems to be the key in the dynamic depletion of the 
vortex stretching for the 2D QG equation. Due to this local geometric
regularity of level sets, Theorem \ref{thm:double} can be applied to 
the 2D QG flow, which explains the double exponential growth of 
$\left| \nabla^{\bot} \theta \right|$ observed in recent numerical 
simulations.

\subsection{Numerical method}

Our numerical simulations are performed using  the pseudo-spectral method 
with the $2/3$ de-aliasing rule in space, and the 4th order classical 
Runge-Kutta in time. We use up to $2048 \times 2048$ space resolution to
resolve the rapidly increasing gradient of $\theta$. The size of the 
time step is determined by the CFL condition. To 
make the time marching more stable, we use $1/6$ of the maximum
allowed CFL number. The computation is done on a 4-CPU (Intel(R) 
Xeon 3.00 GHz) machine with 2048 Kb cache and 6G memory. The
FFT code is from FFTW 3.1.

We use the same initial condition used by Constantin-Majda-Tabak
{\cite{constantin-majda-tabak:1994}}, Ohkitani-Yamada
{\cite{ohkitani-yamada:1997}}, Constantin-Nie-Schorghofer
{\cite{constantin-nie-schorghofer:1998}},
{\cite{constantin-nie-schorghofer:1999}} which contains hyperbolic saddles:
\begin{equation}
  \label{eq:theta.0} \theta_0 \left( x, y \right) = \sin x \sin y + \cos y.
\end{equation}
In Constantin-Nie-Schorghofer {\cite{constantin-nie-schorghofer:1999}}, other
initial conditions were also considered. It was found that the solutions 
corresponding to these initial conditions behave essentially the same as
the one described above. Therefore we decide to focus on the above initial 
condition (\ref{eq:theta.0}) and try to resolve the hyperbolic saddle and the 
fine structure of the direction field of $\nabla^{\bot} \theta$.

We plot the maximum $\left| \nabla^{\bot} \theta \right|$ versus 
time in Figure \ref{fig:loglog}, which is almost identical to Fig. 2 
in Ohkitani-Yamada {\cite{ohkitani-yamada:1997}}. We also compare the 
level set contours obtained by our computations with those obtained in
Constantin-Nie-Schorghofer {\cite{constantin-nie-schorghofer:1998}},
{\cite{constantin-nie-schorghofer:1999}}. They are essentially
indistinguishable.

\begin{center}
  \tmfloat{h}{small}{figure}{\resizebox{3.5in}{!}{\epsfig{file=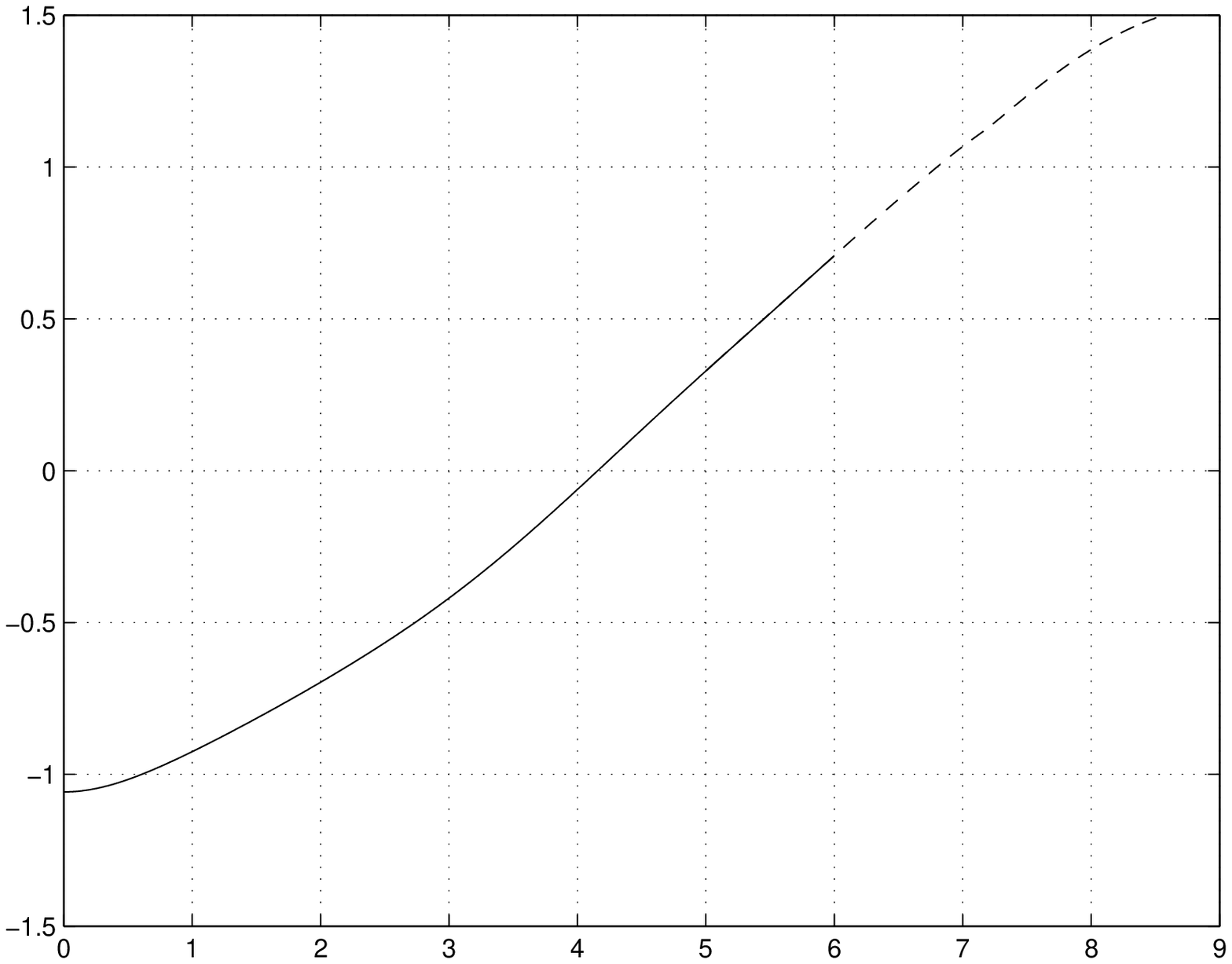}}}{\label{fig:loglog}
  $\tmop{loglog} \left| \nabla^{\bot} \theta \right|$ versus time. The solid
  curve is obtained using $1024 \times 1024$ resolution ($t = 0$ to $t = 6$),
  and the dashed one is obtained using $2048 \times 2048$ resolution (from $t
  = 5$ on). We see that the two curves are almost identical in the time
  interval $\left[ 5, 6 \right]$. }
\end{center}

\subsection{Dynamics of level set geometries}

Next we study the evolution of the two regions mentioned above. We plot, 
at times $t = 5.0,
\hspace{1em} 6.0, \hspace{1em} 6.5, \hspace{1em}$and $7.0$ the 
boundaries of the following two regions:
\begin{enumeratenumeric}
  \item $A_t \equiv \left\{ \left( x, y \right) \mid \left| \nabla^{\bot}
  \theta \left( x, y, t \right) \right| \geqslant \frac{1}{2} \left\|
  \nabla^{\bot} \theta \left( \cdot, t \right) \right\|_{L^{\infty}}
  \right\}$, and
  
  \item $B_t \equiv \left\{ \left( x, y \right) \mid \left| \nabla \xi \right|
  \geqslant 10 \right\}$.
\end{enumeratenumeric}
In Figures \ref{fig:ABt=5.0}, \ref{fig:ABt=6.0}, \ref{fig:ABt=6.5},
\ref{fig:ABt=7.0}, we see that $A_t$ and $B_t$  are essentially disjoint,
although some interlacing of the boundaries can be observed. 
This implies that $\left| \nabla \xi \right|$ 
(consequently $\kappa$ and $\nabla \cdot \xi$) is bounded in the region  
where $\left| \nabla^{\bot} \theta \right|$ achieves its maximum.
This is an interesting result by itself. It says that the local geometric 
property of the level sets is regular in the region of maximum stretching.
Furthermore, it can be seen from Figures \ref{fig:At=5.0}, \ref{fig:At=6.0}, 
\ref{fig:At=6.5}, \ref{fig:At=7.0} that although $A_t$ is severely stretched 
as time increases, the stretching seems to align with the level set curves. 
As a result, we can always pick a level set segment of $O \left( 1 \right)$ 
length, along which the maximum $\left|\nabla^\bot\theta\right|$ is comparable to its 
global maximum. Thus, the conditions of Theorem \ref{thm:double} are satisfied, 
and we can apply Theorem \ref{thm:double} to conclude that the maximum 
growth rate of $\left| \nabla^{\bot} \theta \right|$ is bounded by double 
exponential, which is consistent with recent numerical simulations 
(Ohkitani-Yamada {\cite{ohkitani-yamada:1997}},
Constantin-Nie-Schorghofer {\cite{constantin-nie-schorghofer:1998}},
{\cite{constantin-nie-schorghofer:1999}}).

\begin{center}
  \tmfloat{h}{small}{figure}{\epsfig{file=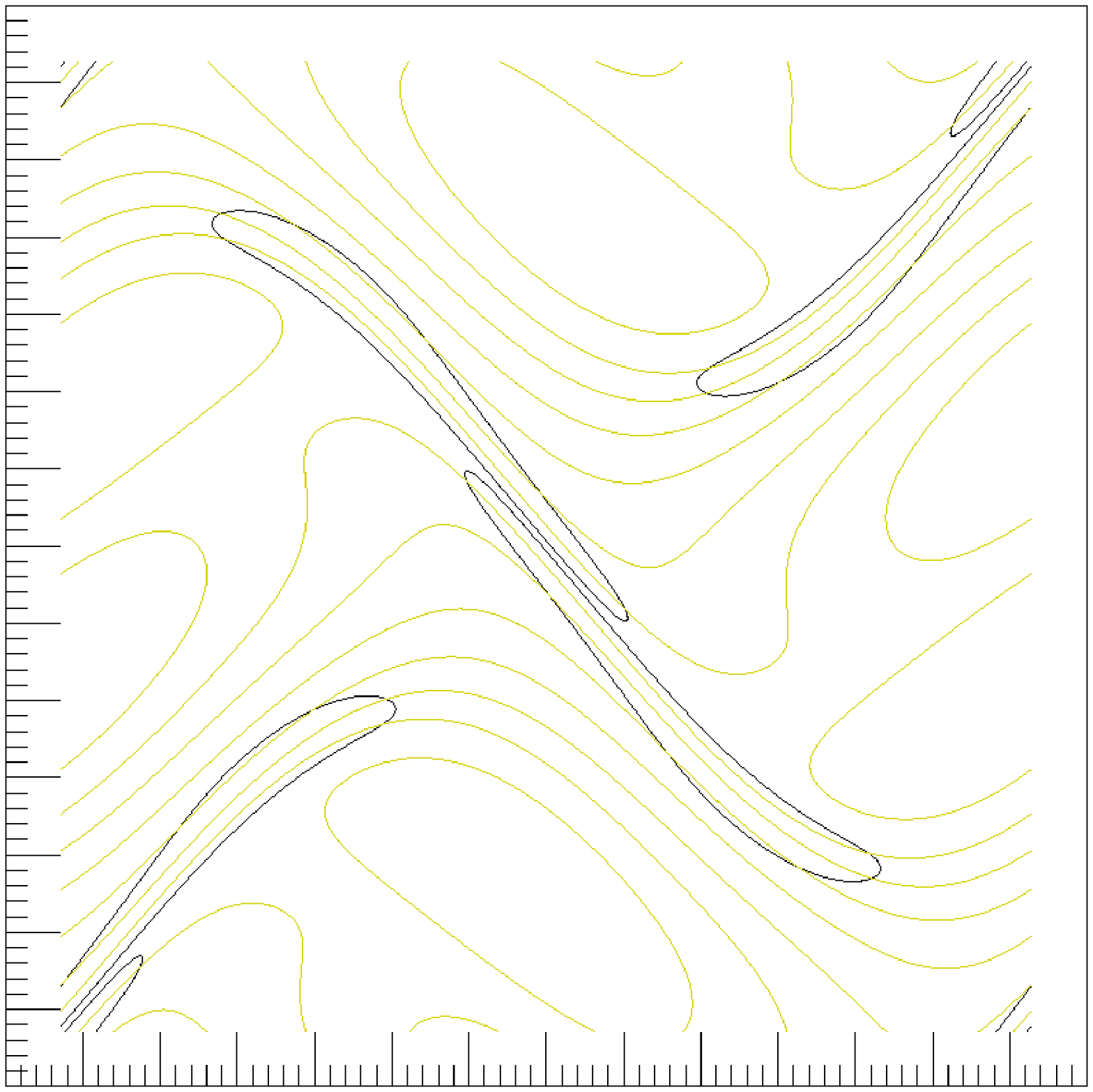}}{\label{fig:At=5.0}
  $A_t$ against the level sets of $\theta$, $t = 5.0$.}
\end{center}

\begin{center}
  \tmfloat{h}{small}{figure}{\epsfig{file=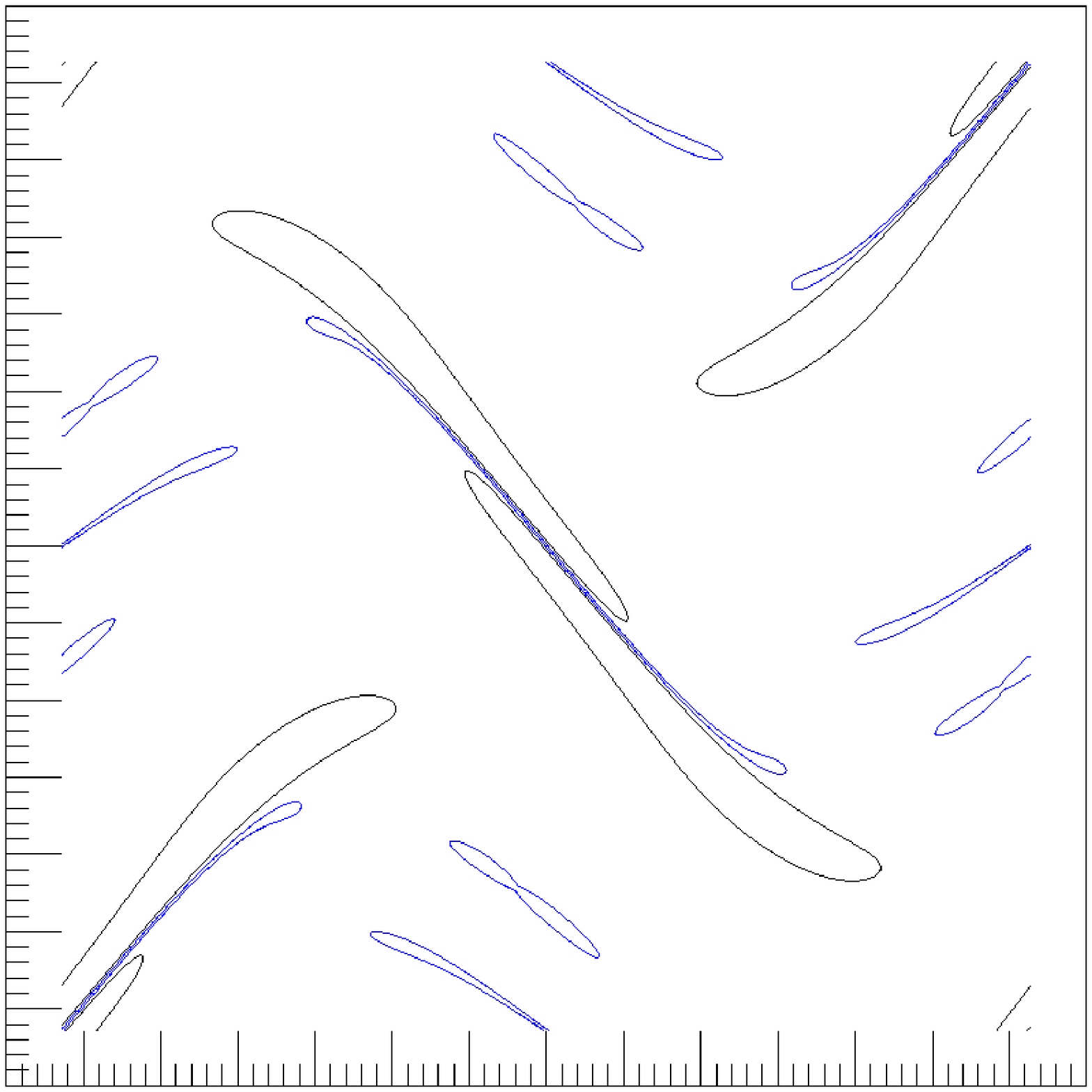}}{\label{fig:ABt=5.0}
  $A_t$ and $B_t$, $t = 5.0$. The boundary of $B_t$ is plotted in blue.}
\end{center}

\begin{center}
  \tmfloat{h}{small}{figure}{\epsfig{file=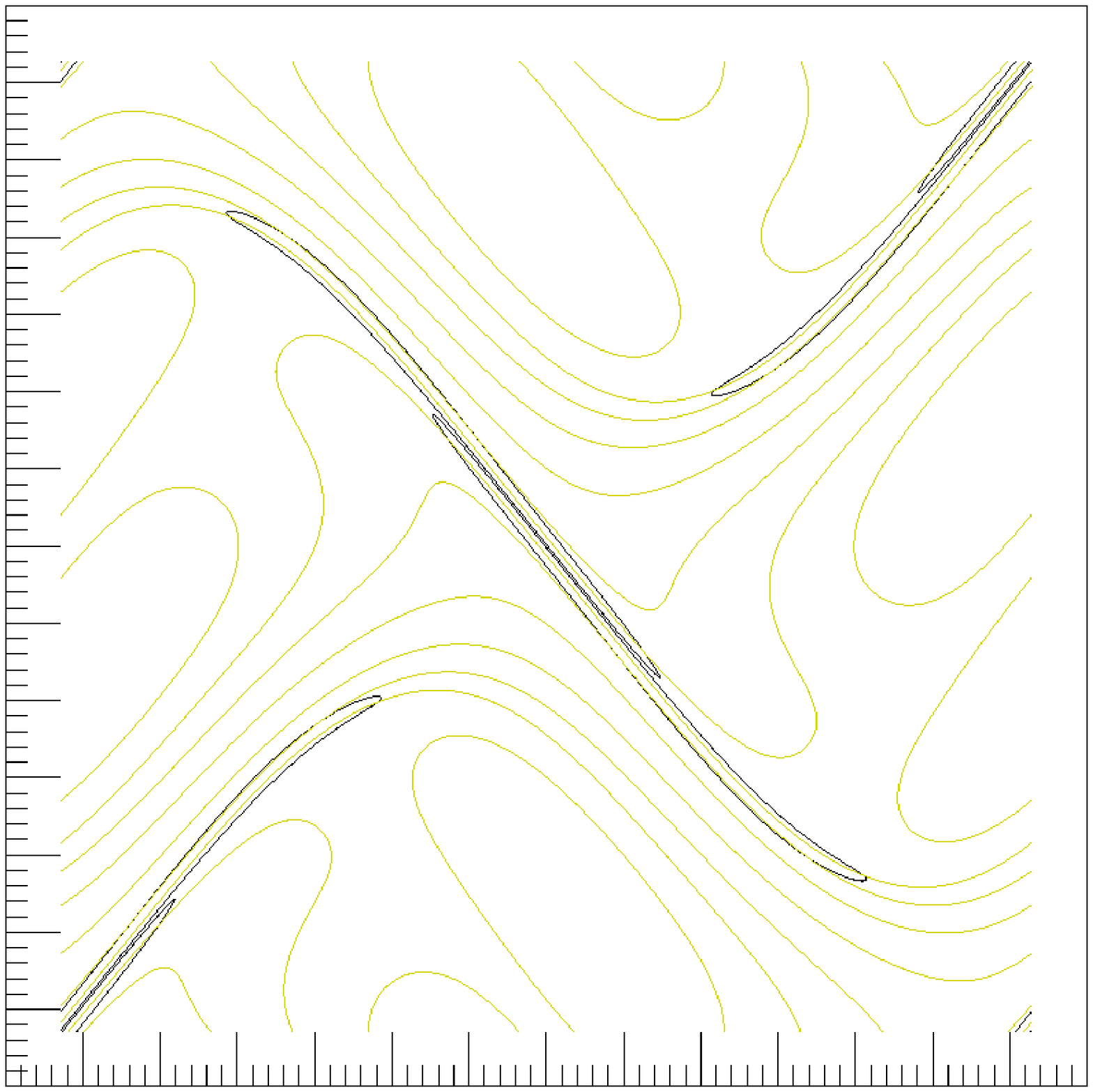}}{\label{fig:At=6.0}
  $A_t$ against level sets of $\theta$, $t = 6.0$}
\end{center}

\begin{center}
  \tmfloat{h}{small}{figure}{\epsfig{file=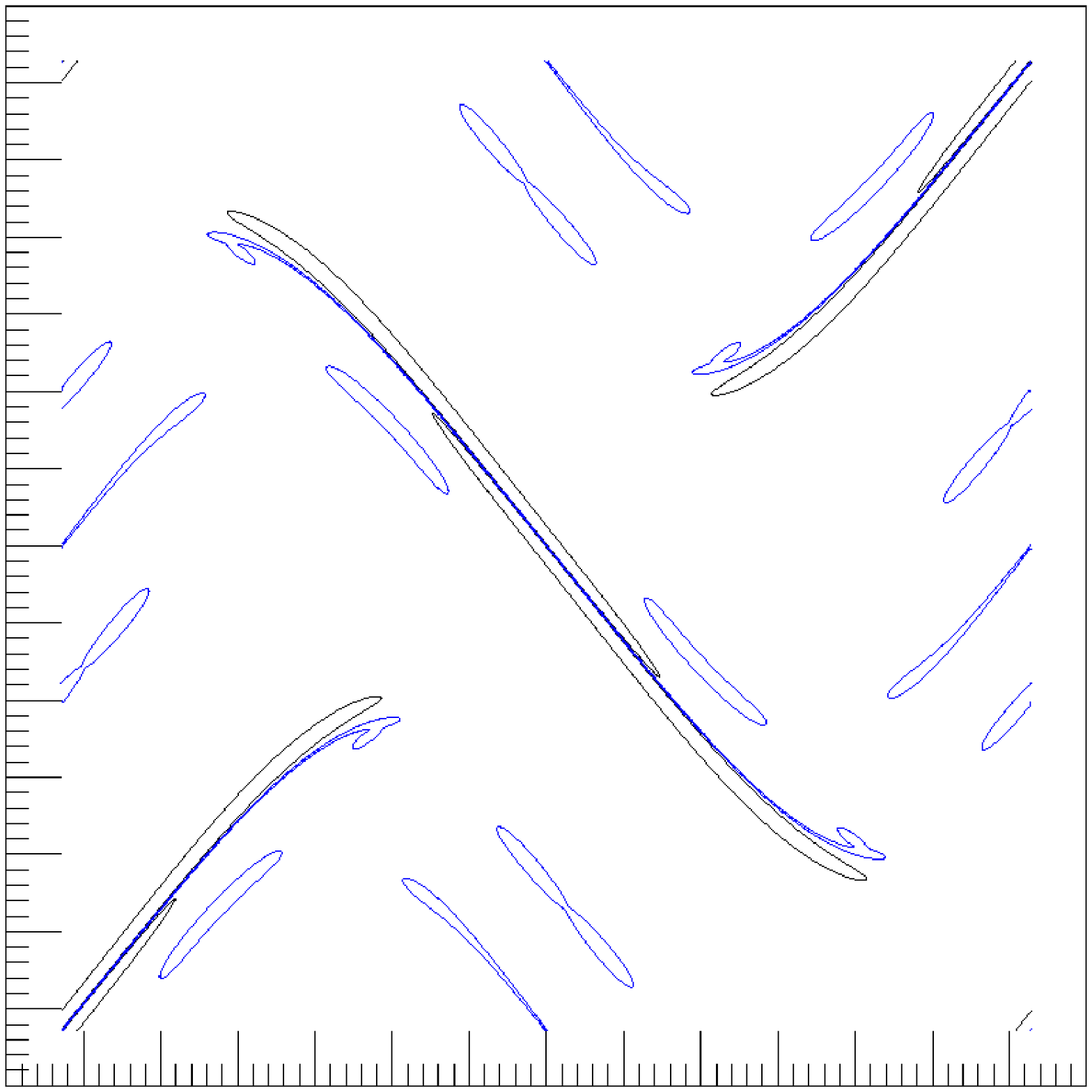}}{\label{fig:ABt=6.0}
  $A_t$ and $B_t$, $t = 6.0$. The boundary of $B_t$ is plotted in blue.}
\end{center}

\begin{center}
  \tmfloat{h}{small}{figure}{\epsfig{file=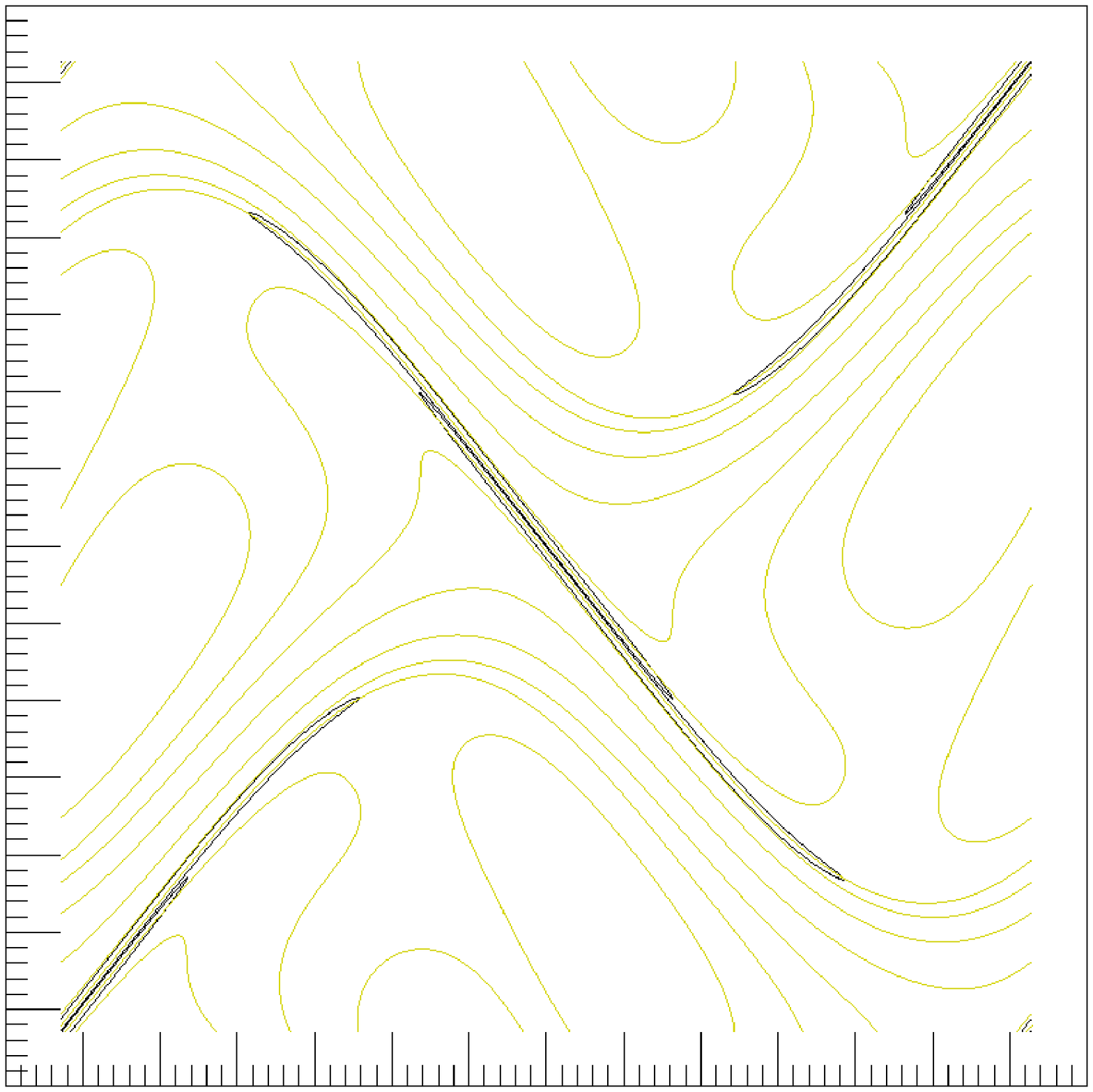}}{\label{fig:At=6.5}
  $A_t$ against level sets of $\theta$, $t = 6.5$}
\end{center}

\begin{center}
  \tmfloat{h}{small}{figure}{\epsfig{file=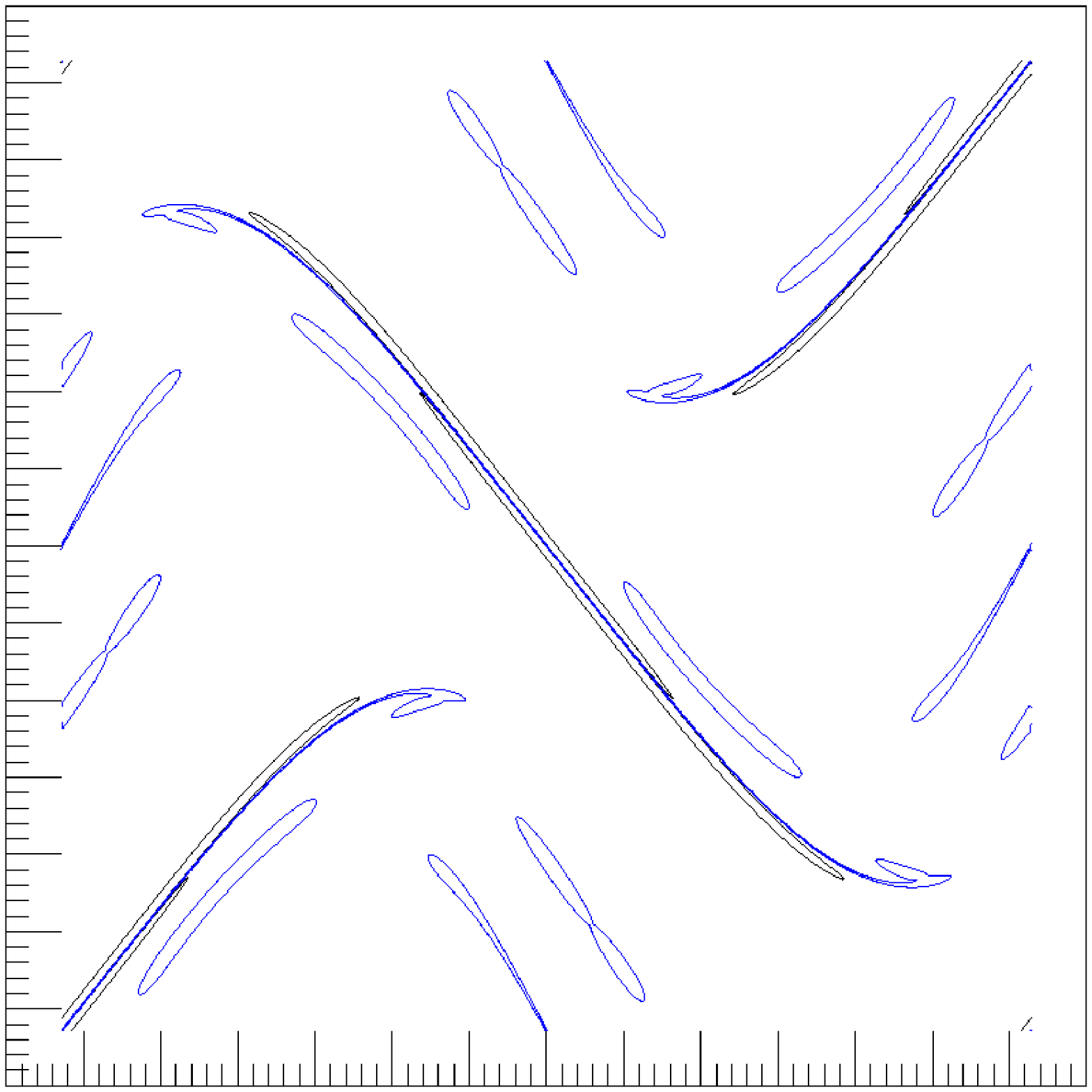}}{\label{fig:ABt=6.5}
  $A_t$ and $B_t$, $t = 6.5$. The boundary of $B_t$ is plotted in blue.}
\end{center}

\begin{center}
  \tmfloat{h}{small}{figure}{\epsfig{file=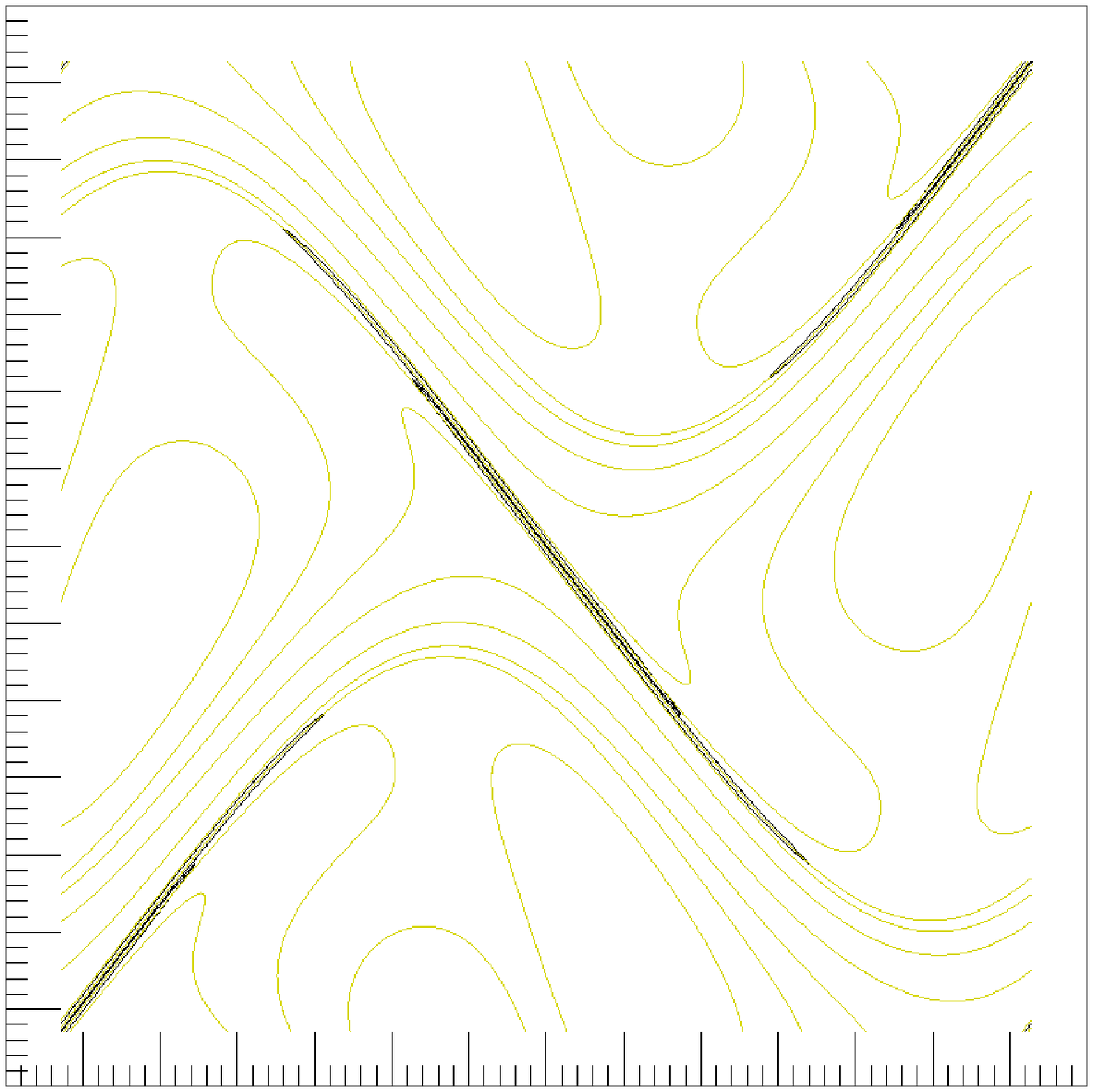}}{\label{fig:At=7.0}
  $A_t$ against level sets of $\theta$, $t = 7.0$}
\end{center}

\begin{center}
  \tmfloat{h}{small}{figure}{\epsfig{file=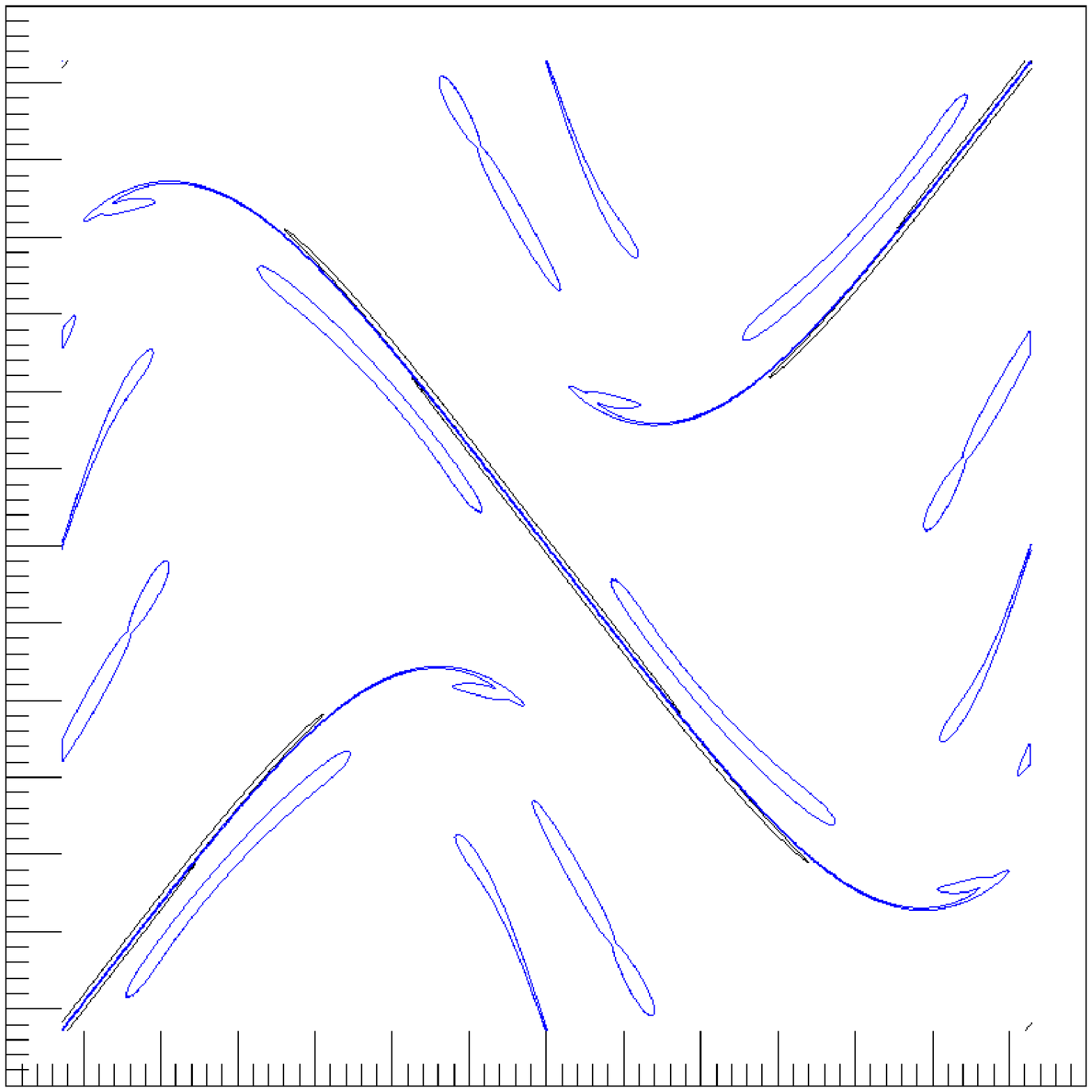}}{\label{fig:ABt=7.0}
  $A_t$ and $B_t$, $t = 7.0$. The boundary of $B_t$ is plotted in blue.}
\end{center}

\begin{center}
  \tmfloat{h}{small}{figure}{\epsfig{file=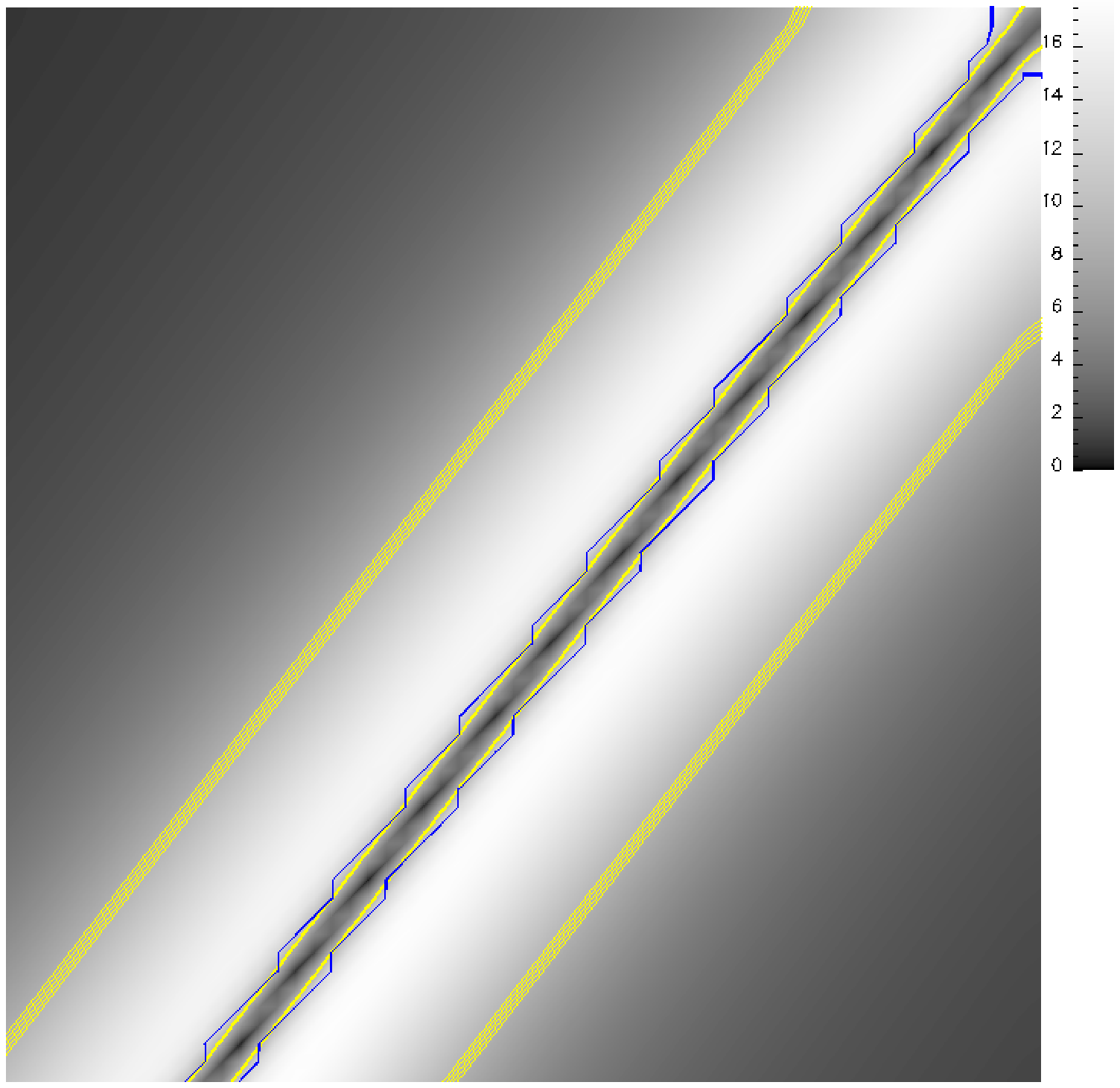}}{\label{fig:ABt=7.0.zoom}Zoom
  in of the upper-right corner of Figure \ref{fig:ABt=7.0} ($t=7.0$). Lighter color corresponds to larger $\left|\nabla^\bot\theta\right|$.
  The boundary of $A_t$ is plotted in yellow (thickened), and the boundary of $B_t$ in blue. We see that
  basically the two sets $A_t$ and $B_t$ miss each other. Note that due to
  periodicity of the data, this figure also reveals what happens in the center part
  of Figure \ref{fig:ABt=7.0}.}
\end{center}

\section{\label{sec:key-est}Level Set Dynamics and Proofs of the Theorems}

\subsection{Level set dynamics and key estimate}

Before proving the theorems, we need to do some preparations. First we fix
notations.
\begin{itemize}
  \item $C$ or $c$: generic constants, whose value may change from line to
  line.
  
  \item $\xi$: the direction of $\nabla^{\bot} \theta \equiv \left( -
  \frac{\partial \theta}{\partial x_2}, \frac{\partial \theta}{\partial x_1}
  \right)^T$, that is, $\xi \equiv \frac{\nabla^{\bot} \theta}{\left|
  \nabla^{\bot} \theta \right|}$ when $\nabla^{\bot} \theta$ does not vanish.
  In the following we will consider level set segments along which $\left|
  \nabla^{\bot} \theta \right|$ are comparable to the maximum $\left\|
  \nabla^{\bot} \theta \right\|_{L^{\infty}}$, therefore we do not need to
  consider the case when $\xi$ is not well-defined.
  
  \item $x, \alpha$: Cartesian coordinate variables. Thus $x, \alpha \in
  \mathbbm{R}^2$.
  
  \item $s, \beta$: arc length variables along the level set under
  consideration.
  
  \item $X \left( \alpha, \tau, t \right)$: the particle trajectory passing
  $\alpha$ at time $\tau$. In other words, $X \left( \alpha, \tau, t \right)$
  solves
  \begin{eqnarray*}
    \frac{\partial X \left( \alpha, \tau, t \right)}{\partial t} & = & u
    \left( X \left( \alpha, \tau, t \right), t \right)\\
    X \left( \alpha, \tau, \tau \right) & = & \alpha .
  \end{eqnarray*}
  For any set $A \subseteq \mathbbm{R}^2$, we denote
  \[ X \left( A, \tau, t \right) \equiv \cup_{\alpha \in A} X \left( \alpha,
     \tau, t \right) . \]
  When $\tau = 0$, we use the traditional notation $X \left( \alpha, t \right)
  \equiv X \left( \alpha, 0, t \right)$.
  
  \item $\sim$: We write $a \left( t \right) \sim b \left( t \right)$ if
  there are absolute constants $c, C > 0$ such that
  \[ c \left| a \left( t \right) \right| \leqslant \left| b \left( t \right)
     \right| \leqslant C \left| a \left( t \right) \right| . \]
  \item $\gtrsim, \lesssim$: We write $a \left( t \right) \gtrsim b \left( t
  \right)$ if there is an absolute constant $c > 0$ such that
  \[ \left| a \left( t \right) \right| \geqslant c \left| b \left( t \right)
     \right| . \]
  $a \left( t \right) \lesssim b \left( t \right)$ is defined similarly. 
\end{itemize}
In the following subsection, we follow the same line of derivation as in
in Deng-Hou-Yu {\cite{deng-hou-yu:2005}}. We  give the details of the 
derivation here to make the presentation self-contained.

\subsubsection{Level Set Dynamics}

First we derive estimates for the stretching of level sets.

\begin{lemma}
  \label{lem:theta.nabla.flow.map}If $\theta$ solves the 2D QG equation with
  initial value $\theta_0$, and if furthermore $X \left( \alpha, t \right)$ is
  the flow map, then we have
  \begin{equation}
    \nabla^{\bot} \theta \left( X \left( \alpha, t \right), t \right) =
    \nabla_{\alpha} X \left( \alpha, t \right) \cdot \nabla_{\alpha}^{\bot}
    \theta_0 \left( \alpha \right) .
  \end{equation}
where the subscript $\alpha$ denotes partial derivative with respect to 
$\alpha$.
\end{lemma}

\begin{proof}
  We prove by a direct calculation.
  \begin{eqnarray*}
    \nabla^{\bot}_{\alpha} \theta_0 \left( \alpha \right) & = &
    \nabla_{\alpha}^{\bot} \theta \left( X \left( \alpha, t \right), t
    \right)\\
    & = & \left( \nabla_{\alpha} \theta \left( X \left( \alpha, t \right), t
    \right) \right)^{\bot}\\
    & = & \left( \nabla_{\alpha} X \cdot \nabla \theta \right)^{\bot}\\
    & = & \left( \nabla_{\alpha} X \right)^{- 1} \cdot \nabla^{\bot} \theta
    \left( X \left( \alpha, t \right), t \right),
  \end{eqnarray*}
  where in the last equality we have used the incompressibility property of
  the flow, i.e. $\det \nabla_{\alpha} X \equiv 1$.
\end{proof}

\begin{lemma}
  \label{lem:s_beta_nabla_theta}Consider a point $X \left( \alpha, t_0, t
  \right)$ carried by the flow. Let $s$ be the arc length variable along the
  level set passing $X \left( \alpha, t_0, t \right)$ at time $t$, and let
  $\beta$ be the arc length variable of the same level set at time $t_0$. Then
  if $\nabla^{\bot} \theta \left( \alpha, t_0 \right) \neq 0$, we have
  \begin{equation}
    \label{eq:s.beta.nabla.theta} \frac{\partial s}{\partial \beta} \left( X
    \left( \alpha, t_0, t \right), t \right) = \frac{\left| \nabla^{\bot}
    \theta \left( X \left( \alpha, t_0, t \right), t \right) \right|}{\left|
    \nabla_{\alpha}^{\bot} \theta \left( \alpha, t_0 \right) \right|}
  \end{equation}
\end{lemma}

\begin{proof}
  By Lemma \ref{lem:theta.nabla.flow.map} we have
  \[ \nabla^{\bot} \theta \left( X \left( \alpha, t \right), t \right) =
     \nabla_{\alpha} X \left( \alpha, t \right) \cdot \nabla_{\alpha}^{\bot}
     \theta_0 \left( \alpha \right) . \]
  Therefore
  \begin{eqnarray*}
    \left| \nabla^{\bot} \theta \left( X, t \right) \right| & = & \xi \left(
    X, t \right) \cdot \nabla^{\bot} \theta \left( X, t \right)\\
    & = & \xi \left( X, t \right) \cdot \nabla_{\alpha} X \left( \alpha, t_0,
    t \right) \cdot \xi \left( \alpha, t_0 \right)  \left|
    \nabla_{\alpha}^{\bot} \theta \left( \alpha, t_0 \right) \right|\\
    & = & \xi \left( X, t \right) \cdot \left[ \xi \left( \alpha, t_0 \right)
    \cdot \nabla_{\alpha} X \left( \alpha, t_0, t \right) \right]  \left|
    \nabla_{\alpha}^{\bot} \theta \left( \alpha, t_0 \right) \right|\\
    & = & \xi \left( X, t \right) \cdot \frac{\partial X}{\partial \beta} 
    \left| \nabla_{\alpha}^{\bot} \theta \left( \alpha, t_0 \right) \right|\\
    & = & \xi \left( X, t \right) \cdot \frac{\partial X}{\partial s}  \left|
    \nabla_{\alpha}^{\bot} \theta \left( \alpha, t_0 \right) \right| 
    \frac{\partial s}{\partial \beta}\\
    & = & \left[ \xi \left( X, t \right) \cdot \xi \left( X, t \right)
    \right]  \left| \nabla_{\alpha}^{\bot} \theta \left( \alpha, t_0 \right)
    \right|  \frac{\partial s}{\partial \beta}\\
    & = & \left| \nabla_{\alpha}^{\bot} \theta \left( \alpha, t_0 \right)
    \right|  \frac{\partial s}{\partial \beta},
  \end{eqnarray*}
  where we have used the fact $\xi \left( X, t \right) = \frac{\partial X
  \left( \alpha, t_0, t \right)}{\partial s}$ where $s$ is the arc length
  variable at time $t$, and in particular, $\xi \left( \alpha, t_0 \right) =
  \frac{\partial \alpha}{\partial \beta}$ where $\beta$ is the arc length
  variable at time $t_0$. 
\end{proof}

Now we can write the evolution equation for $\frac{\partial s}{\partial
\beta}$ (in the following denoted as $s_{\beta}$). In Constantin-Majda-Tabak
{\cite{constantin-majda-tabak:1994}}, it is derived that
\begin{equation}
  D_t \left| \nabla^{\bot} \theta \right| = \left( \xi \cdot \nabla u \cdot
  \xi \right)  \left| \nabla^{\bot} \theta \right| .
\end{equation}
Thanks to (\ref{eq:s.beta.nabla.theta}) we immediately have the evolution
equation for $s_{\beta}$:
\begin{equation}
  D_t \left( s_{\beta} \right) = \left( \xi \cdot \nabla u \cdot \xi \right)
  s_{\beta} .
\end{equation}
Now observe
\begin{eqnarray*}
  \xi \cdot \nabla u \cdot \xi & = & \left( \xi \cdot \nabla \right) \left( u
  \cdot \xi \right) - u \cdot \left( \xi \cdot \nabla \right) \xi\\
  & = & \left( u \cdot \xi \right)_s - \kappa \left( u \cdot \tmmathbf{n}
  \right)
\end{eqnarray*}
where we have used $\xi \cdot \nabla = \partial / \partial s$ and the 
Frenet relation
\[ \frac{\partial \xi}{\partial s} = \kappa \tmmathbf{n} \]
with $\kappa = \left| \xi \cdot \nabla \xi \right|$ being the curvature, and
$\tmmathbf{n}$ the unit normal vector of the level set curve. Therefore we
have an alternative formulation of the evolution equation of $s_{\beta}$:
\begin{equation}
  \label{eq:s.beta.evol} D_t \left( s_{\beta} \right) = \left( u \cdot \xi
  \right)_{\beta} - \kappa \left( u \cdot \tmmathbf{n} \right) s_{\beta} .
\end{equation}

Next we consider a small level set segment at time $t_0$, whose two ends are
denoted by arc lengths $\beta_1 < \beta_2$. Let $l_t$ denote this level set
segment at time $t$, and let $l \left( t \right)$ denote its length at time
$t$. To study its stretching over time, we integrate (\ref{eq:s.beta.evol})
along $l_{t_0}$ to obtain
\begin{eqnarray*}
  D_t \left[ s \left( \beta_2, t \right) - s \left( \beta_1, t \right) \right]
  & = & \left( u \cdot \xi \right) \left( X \left( \beta_2, t_0, t \right), t
  \right) - \left( u \cdot \xi \right) \left( X \left( \beta_1, t_0, t
  \right), t \right)\\
  &  & - \int_{\beta_1}^{\beta_2} \left[ \kappa \left( u \cdot \tmmathbf{n}
  \right) s_{\beta} \right] \left( X \left( \beta, t_0, t \right), t \right)
  \mathd \beta \hspace*{\fill}   \text{\tmtextup{(1)}}.
\end{eqnarray*}
Now we further integrate from $t_0$ to some later time $t$. We get
\begin{eqnarray*}
  s \left( \beta_2, t \right) - s \left( \beta_1, t \right) & = & s \left(
  \beta_2, t_0 \right) - s \left( \beta_1, t_0 \right)\\
  &  & + \int_{t_0}^t \left[ \left( u \cdot \xi \right) \left( X \left(
  \beta_2, t_0, \tau \right), \tau \right) - \left( u \cdot \xi \right) \left(
  X \left( \beta_1, t_0, \tau \right), \tau \right) \right] \mathd \tau\\
  &  & - \int_{t_0}^t \int_{\beta_1}^{\beta_2} \left[ \kappa \left( u \cdot
  \tmmathbf{n} \right) s_{\beta} \right] \left( X \left( \beta, t_0, \tau
  \right), \tau \right) \mathd \beta \mathd t \hspace*{\fill}  
  \text{\tmtextup{(2)}}.
\end{eqnarray*}
Using the notation $l \left( t \right) \equiv s \left( \beta_2, t \right) - s
\left( \beta_1, t \right)$, we obtain
\begin{equation}
  \label{eq:l.stretching} l \left( t \right) \leqslant l \left( t_0 \right) +
  \int_{t_0}^t \left[ u_{\xi} \left( \tau \right) + k \left( \tau \right) u_n
  \left( \tau \right) l \left( \tau \right) \right] \mathd \tau
\end{equation}
where $k \left( \tau \right) = \max_{l_t} \kappa$, $u_{\xi} \left( \tau
\right) = \max_{x, y \in l_{\tau}} \left| \left( u \cdot \xi \right) \left( x,
\tau \right) - \left( u \cdot \xi \right) \left( y, \tau \right) \right|$ and
$u_n \left( \tau \right) = \max_{x \in l_{\tau}} \left| \left( u \cdot
\tmmathbf{n} \right) \left( x, \tau \right) \right|$. Denoting $U \left( \tau
\right) = \left\| u \left( \cdot, \tau \right) \right\|_{L^{\infty} \left(
\mathbbm{R}^2 \right)}$, we get the following weaker estimate which is enough
for our purpose in this paper:
\begin{equation}
  \label{eq:l.stretching.weak} l \left( t \right) \leqslant l \left( t_0
  \right) + 2 \int_{t_0}^t \left[ 1 + k \left( \tau \right) l \left( \tau
  \right) \right] U \left( \tau \right) \mathd \tau .
\end{equation}

\subsubsection{Estimate of $\left| \nabla^{\bot} \theta \right|$ growth}

To apply (\ref{eq:l.stretching.weak}) to the estimate of $\left| \nabla^{\bot}
\theta \right|$, we need to relate the stretching of $l_t$ to the growth of
$\left| \nabla^{\bot} \theta \right|$. This is given by the following two
lemmas.

\begin{lemma}
  \label{lem:div-xi-nabla-theta}Let $\xi \left( x, t \right)$ be the direction
  of $\nabla^{\bot} \theta$. Assume at some time $t$, the solution $\theta
  \left( x, t \right)$ is $C^1$ in $x$. Then at this time $t$, for any $x$
  such that $\nabla \theta \neq 0$, there holds
  \begin{equation}
    \label{eq:div.free.implication} \frac{\partial \left| \nabla^{\bot} \theta
    \right|}{\partial s} \left( x, t \right) = - \left( \left( \nabla \cdot
    \xi \right)  \left| \nabla^{\bot} \theta \right| \right) \left( x, t
    \right),
  \end{equation}
  where $s$ is the arc length variable along the level set passing $x$ at time
  $t$.
  
  Furthermore, if we denote this vortex line by $l$, then for any $y \in l$
  such that $\nabla^{\bot} \theta$ does not vanish at any point in the level
  set segment connecting $x$ and $y$, we have
  \begin{equation}
    \label{eq:nabla.theta.x.y} \left| \nabla^{\bot} \theta \right| \left( y, t
    \right) = \left| \nabla^{\bot} \theta \right| \left( x, t \right)
    e^{\int_x^y \left( - \nabla \cdot \xi \right) \mathd s},
  \end{equation}
  where the integration is along $l$. 
\end{lemma}

\begin{proof}
  We compute
  \begin{eqnarray*}
    0 = \nabla \cdot \left( \nabla^{\bot} \theta \right) & = & \nabla \cdot
    \left( \left| \nabla^{\bot} \theta \right| \xi \right)\\
    & = & \left( \xi \cdot \nabla \right) \left| \nabla^{\bot} \theta \right|
    + \left( \nabla \cdot \xi \right)  \left| \nabla^{\bot} \theta \right|
  \end{eqnarray*}
  which immediately gives (\ref{eq:div.free.implication}) after noticing $\xi
  \cdot \nabla = \frac{\partial}{\partial s}$. (\ref{eq:nabla.theta.x.y}) is
  the immediate result of solving (\ref{eq:div.free.implication}) along $l$. 
\end{proof}

\begin{lemma}
  \label{lem:Omega.l.relation}Let $l_t$ be a level set segment carried by the
  flow, i.e., $l_t = X \left( l_{t_0}, t_0, t \right)$ for some earlier time
  $t_0$. Define
  \[ m \left( t \right) \equiv \max_{x \in l_t} \left| \left( \nabla \cdot
     \xi \right) \left( x, t \right) \right|, \]
  where $\xi = \nabla^{\bot} \theta / \left| \nabla^{\bot} \theta \right|$ is
  the unit tangent vector. If we further denote $\Omega_l \left( t \right)
  \equiv \max_{l_t} \left| \nabla^{\bot} \theta \right|$, then the following
  inequalities hold:
  \begin{equation}
    \label{eq:Omega.l.relation} e^{- m \left( t \right) l \left( t \right)}
    \frac{\Omega_l \left( t \right)}{\Omega_l \left( t_0 \right)} \leqslant
    \frac{l \left( t \right)}{l \left( t_0 \right)} \leqslant e^{m \left( t_0
    \right) l \left( t_0 \right)} \frac{\Omega_l \left( t \right)}{\Omega_l
    \left( t_0 \right)} .
  \end{equation}
\end{lemma}

\begin{proof}
  Recall that $\beta$ is the arc length variable at time $t_0$. We have
  \begin{eqnarray*}
    l \left( t \right) & = & \int_{\beta_1}^{\beta_2} s_{\beta} \mathd \beta\\
    & = & \int_{\beta_1}^{\beta_2} \frac{\left| \nabla^{\bot} \theta \left(
    \alpha, t_0, t \right) \right|}{\left| \nabla^{\bot} \theta \left( \alpha,
    t_0 \right) \right|} \mathd \beta\\
    & \leqslant & \int_{\beta_1}^{\beta_2} \frac{\Omega_l \left( t
    \right)}{e^{- m \left( t_0 \right) l \left( t_0 \right)} \Omega_l \left(
    t_0 \right)} \mathd \beta\\
    & = & e^{m \left( t_0 \right) l \left( t_0 \right)} \frac{\Omega_l \left(
    t \right)}{\Omega_l \left( t_0 \right)} l \left( t_0 \right),
  \end{eqnarray*}
  where the inequality is a direct result of (\ref{eq:nabla.theta.x.y}).
  
  The other inequality is proved similarly. 
\end{proof}

Finally, combining (\ref{eq:l.stretching.weak}) and
(\ref{eq:Omega.l.relation}), we obtain the following estimate (after dividing
both sides by $l \left( t_0 \right)$).
\begin{equation}
  \label{eq:key} \Omega_l \left( t \right) \leqslant e^{m \left( t \right) l
  \left( t \right)} \Omega_l \left( t_0 \right)  \left[ 1 + \frac{2}{l \left(
  t_0 \right)}  \int_{t_0}^t \left[ 1 + k \left( \tau \right) l \left( \tau
  \right) \right] U \left( \tau \right) \mathd \tau \right] ,
\end{equation}
where  $\Omega_l \left( t \right) \equiv \max_{l_t} \left| \nabla^{\bot}
\theta \right|$, $m \left( t \right) \equiv \max_{x \in l_t} \left| \left(
\nabla \cdot \xi \right) \left( x, t \right) \right|$, $k \left( \tau \right)
= \max_{l_t} \kappa$, and $U \left( \tau \right) = \left\| u \left( \cdot,
\tau \right) \right\|_{L^{\infty} \left( \mathbbm{R}^2 \right)}$. This
estimate will play a key role in the proofs of the theorems.

\subsection{Proofs of the Theorems}

\subsubsection{Proof of Theorem \ref{thm:triple}. }

From estimate (\ref{eq:key}), it is clear that an estimate of $U \left( \tau
\right)$ is needed. Such an estimate is derived in Cordoba
{\cite{cordoba:1998}}. We summarize his estimate into the following lemma.

\begin{lemma}
  \label{lem:cordoba}There exists a generic constant $C > 0$ such that for $t
  > 0$
  \[ \left\| u \left( \cdot, t \right) \right\|_{L^{\infty}} \leqslant C \log
     \left\| \nabla^{\bot} \theta \left( \cdot, t \right)
     \right\|_{L^{\infty}} \]
  provided that $\left\| \nabla^{\bot} \theta \left( \cdot, t \right)
  \right\|_{L^{\infty}} > e$.
\end{lemma}

In our notations, the above estimate is just
\begin{equation}
  \label{eq:U.Omega.est} U \left( t \right) \leqslant C \log \Omega \left( t
  \right).
\end{equation}
Also, we do not need to worry about the condition $\Omega \left( t \right) >
e$ since we will consider level set segments carrying large $\left|
\nabla^{\bot} \theta \right|$.

Now we are ready to prove Theorem \ref{thm:triple}. First we give a heuristic
``proof''. By the assumptions of the theorem, we have $m \left( t \right) l
\left( t \right) \leqslant C_0$, and $\Omega_l \left( t \right) \geqslant c_0
\Omega \left( t \right)$. Thus letting $R = e^{C_0} / c_0$, we have
\begin{equation}
  \label{eq:key.R} \Omega \left( t \right) \leqslant R \Omega \left( t_0
  \right)  \left[ 1 + \frac{C}{l \left( t_0 \right)} \int_{t_0}^t \left[ \log
  \Omega \left( \tau \right) + 1 \right] \mathd \tau \right]
\end{equation}
where we have used (\ref{eq:U.Omega.est}). Intuitively, after taking one
derivative with respect to $t$, and then setting $t_0 = t$, we would get
\[ \Omega' \left( t \right) \leqslant C \Omega \left( t \right) \log \Omega
   \left( t \right) \tmop{loglog} \Omega \left( t \right). \]
This would give the triple exponential bound. Note that in getting the above
differential inequality, we have naively estimated $l \left( t_0 \right)$ by
$\frac{c_L}{\tmop{loglog} \Omega \left( t \right)}$. This estimate cannot be
derived directly from the assumption $l \left( t \right) \geqslant
\frac{c_L}{\tmop{loglog} \Omega \left( t \right)}$ since in general we only
have $l \left( t_0 \right) < l \left( t \right)$. Therefore we need to bound
$l \left( t_0 \right)$ from below using $\Omega \left( t \right)$ . In the
following, we will obtain this lower bound (when $t_0$ and $t$ are not far
apart) and establish the triple exponential upper bound rigorously.

First we outline the main steps.
\begin{itemize}
  \item {\tmstrong{Outline of the main steps}}. The proof consists of four
  steps.
  \begin{enumeratenumeric}
    \item Divide $\left[ T_0, T_{\ast} \right)$ into intervals $\left[ t_k,
    t_{k + 1} \right)$ such that
    \begin{equation}
      \label{eq:Omega.ratio.r} \frac{\Omega \left( t_{k + 1} \right)}{\Omega
      \left( t_k \right)} = r
    \end{equation}
    for some constant $r > R \equiv e^{C_0} / c_0$. One of the reasons for
    doing this partition is to obtain an sharp lower bound estimate for $l
    \left( t_k \right)$ within each time interval $\left[ t_k, t_{k + 1}
    \right)$ using our relationship between the relative growth of $\Omega
    \left( t \right)$ and the relative growth of arc length stretching between
    two different times.
    
    \item Use (\ref{eq:Omega.l.relation}) to obtain a lower bound estimate
    for $l \left( t_k \right)$, which in turn gives an upper bound for $\Omega
    \left( t_{k + 1} \right)$:
    \begin{equation}
      \label{eq:Omega.t.kp1.t.k} \Omega \left( t_{k + 1} \right) \leqslant R
      \Omega \left( t_k \right)  \left[ 1 + C \frac{\left( 1 + C_0 \right)
      Rr}{c_L} \tmop{loglog} \Omega \left( t_k \right)  \int_{t_k}^{t_{k + 1}}
      \left[ \log \Omega \left( \tau \right) + 1 \right] \mathd \tau \right] .
    \end{equation}
    \item Use (\ref{eq:Omega.t.kp1.t.k}) to obtain a local triple exponential
    estimate for $\Omega \left( t_{k + 1} \right)$:
    \begin{eqnarray*}
      \tmop{logloglog} \Omega \left( t_{k + 1} \right) & \leqslant &
      \tmop{logloglog} \Omega \left( t_{k + 1} \right) + C \frac{R^2 r \left(
      1 + C_0 \right)}{c_L}  \left( t_{k + 1} - t_k \right)\\
      &  & + \frac{\log R}{\log \Omega \left( t_k \right) \tmop{loglog}
      \Omega \left( t_k \right)} . \hspace*{\fill}   \text{\tmtextup{(3)}}
      \label{eq:triple.est.local}
    \end{eqnarray*}
    \item Sum up the estimates for each $\left[ t_k, t_{k + 1} \right)$ to
    obtain:
    \begin{eqnarray*}
      \tmop{logloglog} \Omega \left( t_n \right) & \leqslant &
      \tmop{logloglog} \Omega \left( t_0 \right) + C \frac{R^2 r \left( 1 +
      C_0 \right)}{c_L}  \left( t_n - t_0 \right)\\
      &  & + \sum_{i = 0}^{n - 1} \frac{\log R}{\log \Omega \left( t_i
      \right) \tmop{loglog} \Omega \left( t_i \right)} . \hspace*{\fill}  
      \text{\tmtextup{(4)}} \label{eq:triple.est.global}
    \end{eqnarray*}
    It can be shown that the sum in the right hand side of
    (\ref{eq:triple.est.global}) can be bounded as follows:
    \begin{equation}
      \label{eq:triple.rhs} \sum_{i = 0}^{n - 1} \frac{\log R}{\log \Omega
      \left( t_i \right) \tmop{loglog} \Omega \left( t_i \right)} \leqslant
      \frac{\log R}{\log r} \tmop{logloglog} \Omega \left( t_n \right) + C
    \end{equation}
    for some constant $C$. Now substituting (\ref{eq:triple.rhs}) into
    (\ref{eq:triple.est.global}) would give the desired triple exponential
    estimate for $\Omega \left( t_n \right)$:
    \begin{equation}
      \tmop{logloglog} \Omega \left( t_n \right) \leqslant \frac{\log r}{\log
      r - \log R}  \left[ C \frac{R^2 r \left( 1 + C_0 \right)}{c_L}  \left(
      t_n - t_0 \right) + C' \right] .
    \end{equation}
  \end{enumeratenumeric}
\end{itemize}
Now we carry out the above four steps in detail.
\begin{itemize}
  \item {\tmstrong{Partition of the time interval.}}
  
  Let $r$ be any constant such that $r > R$ and $t_0 \in \left[ T_0, T_{\ast}
  \right)$ close enough to $T$ so that $\Omega \left( t_0 \right) > 2 e$ and
  $\tmop{loglog} \left( r \Omega \left( t_0 \right) \right) \leqslant 2
  \tmop{loglog} \Omega \left( t_0 \right)$. Define $t_0 < t_1 < \cdots < t_k <
  \cdots < T_{\ast}$ by (\ref{eq:Omega.ratio.r}), i.e.
  \begin{equation}
    \label{eq:Omega.ratio.r.copy} \frac{\Omega \left( t_{k + 1}
    \right)}{\Omega \left( t_k \right)} = r.
  \end{equation}
  If there exists $n \in \mathbbm{N}$ such that we cannot find $t_{n + 1}$
  using (\ref{eq:Omega.ratio.r.copy}), or equivalently, such that for any $t
  \in \left( t_n, T^{\ast} \right)$,
  \[ \frac{\Omega \left( t) \right.}{\Omega \left( t_n \right)} < r, \]
  then $\Omega \left( t \right)$ remains bounded in $\left[ 0, T_{\ast}
  \right]$, and thus no blow-up can occur. Therefore we assume that for all $k
  \in \mathbbm{N}$ we can find $t_k$ iteratively such that
  (\ref{eq:Omega.ratio.r.copy}) is satisfied. Since $\lim_{k \nearrow \infty}
  \Omega \left( t_k \right) = \infty$ and $T_{\ast}$ is the smallest time such
  that $\int_0^{T_{\ast}} \Omega \left( \tau \right) \mathd \tau = \infty$
  according to the BKM type criterion (\ref{eq:qg.BKM}) derived in
  Constantin-Majda-Tabak {\cite{constantin-majda-tabak:1994}}, we must have
  $t_k \nearrow T_{\ast}$.
  
  \item {\tmstrong{Estimate of the lower bound for $l \left( t_k \right)$.}}
  
  We apply (\ref{eq:key.R}) to the time interval $\left[ t_k, t_{k + 1}
  \right]$. for any $t \in \left[ t_k, t_{k + 1} \right]$, choose $l_{t_{k +
  1}} \subset L_{t_{k + 1}}$ so that $\Omega_l \left( t_{k + 1} \right) =
  \Omega_L \left( t_{k + 1} \right)$, and $l \left( t_{k + 1} \right) =
  \frac{c_L}{\tmop{loglog} \Omega \left( t_{k + 1} \right)}$, and let $l_t$ be
  such that $l_{t_{k + 1}} = X \left( l_t, t, t_{k + 1} \right)$, i.e., $l_t$
  is the pullback of $l_{t_{k + 1}}$ to time $t \in \left[ t_k, t_{k + 1}
  \right]$. By the assumptions of Theorem \ref{thm:triple} we have $l_t
  \subset L_t$ for all $t \in \left[ t_k, t_{k + 1} \right]$. Therefore,
  \begin{equation}
    \Omega \left( t \right) \leqslant R \Omega \left( t_k \right)  \left[ 1 +
    C \frac{1 + C_0}{l \left( t_k \right)}  \int_{t_k}^t \left[ \log \Omega
    \left( \tau \right) + 1 \right] \mathd \tau \right] .
  \end{equation}
  Next we obtain a lower bound for $l \left( t_k \right)$. Using
  (\ref{eq:Omega.l.relation}) we have
  \[ \frac{l \left( t_{k + 1} \right)}{l \left( t_k \right)} \leqslant R
     \frac{\Omega \left( t_{k + 1} \right)}{\Omega \left( t_k \right)} = Rr,
  \]
  which gives
  \[ \frac{1}{l \left( t_k \right)} \leqslant \frac{Rr}{l \left( t_{k + 1}
     \right)} = \frac{Rr \tmop{loglog} \Omega \left( t_{k + 1} \right)}{c_L}
     \leqslant \frac{2 Rr \tmop{loglog} \Omega \left( t_k \right)}{c_L} \]
  since $\Omega \left( t_k \right) > \Omega \left( t_0 \right)$ is large
  enough by our choice of $t_0$. Thus we obtain the upper bound
  \begin{equation}
    \label{eq:key.est.tk} \Omega \left( t \right) \leqslant R \Omega \left(
    t_k \right)  \left[ 1 + C \frac{\left( 1 + C_0 \right) Rr}{c_L}
    \tmop{loglog} \Omega \left( t_k \right)  \int_{t_k}^t \left[ \log \Omega
    \left( \tau \right) + 1 \right] \mathd \tau \right] .
  \end{equation}
  for all $t \in \left[ t_k, t_{k + 1} \right]$, where $C$ is some absolute
  constant independent of any parameters.
  
  \item {\tmstrong{Local triple exponential estimate.}}
  
  Define $\tilde{\Omega} \left( t \right)$, $t \in \left[ t_k, t_{k + 1}
  \right]$ by
  \begin{equation}
    \label{eq:def.tilde.Omega} \tilde{\Omega} \left( t \right) = R \Omega
    \left( t_k \right)  \left[ 1 + C \frac{\left( 1 + C_0 \right) Rr}{c_L}
    \tmop{loglog} \Omega \left( t_k \right)  \int_{t_k}^{t_{k + 1}} \left[
    \log \tilde{\Omega} \left( \tau \right) + 1 \right] \mathd \tau \right] .
  \end{equation}
  First we prove $\Omega \left( t \right) < \tilde{\Omega} \left( t \right)$
  for all $t \in \left[ t_k, t_{k + 1} \right]$. When $t = t_k$ we have
  $\tilde{\Omega} \left( t_k \right) = R \Omega \left( t_k \right) > \Omega
  \left( t_k \right)$. Now suppose that there exists $\delta \in \left( 0,
  t_{k + 1} - t_k \right]$ so that $\tilde{\Omega} \left( t \right) > \Omega
  \left( t \right)$ when $t \in \left[ t_k, t_k + \delta \right)$, and
  $\tilde{\Omega} \left( t_k + \delta \right) = \Omega \left( t_k + \delta
  \right)$. Using (\ref{eq:key.est.tk}) and substituting $\tilde{\Omega}
  \left( t_k + \delta \right) = \Omega \left( t_k + \delta \right)$ into
  (\ref{eq:def.tilde.Omega}), we obtain
  \[ \int_{t_k}^{t_k + \delta} \log \tilde{\Omega} \left( \tau \right) \mathd
     \tau \leqslant \int_{t_k}^{t_k + \delta} \log \Omega \left( \tau \right)
     \mathd \tau \]
  which contradicts the assumption that $\tilde{\Omega} \left( t \right) >
  \Omega \left( t \right)$ when $t \in \left[ t_k, t_k + \delta \right)$!
  Therefore, such $\delta$ cannot exist, which means $\Omega \left( t \right)
  < \tilde{\Omega} \left( t \right)$ for all $t \in \left[ t_k, t_{k + 1}
  \right]$.
  
  Next we differentiate (\ref{eq:def.tilde.Omega}) with respect to $t$ and
  get
  \[ \tilde{\Omega}' \left( t \right) = C \frac{R^2 r \left( 1 + C_0
     \right)}{c_L} \Omega \left( t_k \right) \tmop{loglog} \Omega \left( t_k
     \right)  \left[ \log \tilde{\Omega} \left( t \right) + 1 \right] . \]
  Using $\tilde{\Omega} \left( t \right) > \Omega \left( t \right)$, we easily
  obtain
  \begin{eqnarray*}
    \left( \tmop{logloglog} \tilde{\Omega} \left( t \right) \right)' & = & C
    \frac{R^2 r \left( 1 + C_0 \right)}{c_L}  \frac{\Omega \left( t_k \right)
    \tmop{loglog} \Omega \left( t_k \right)  \left[ \log \tilde{\Omega} \left(
    t \right) + 1 \right]}{\tilde{\Omega} \left( t \right) \tmop{loglog}
    \tilde{\Omega} \left( t \right) \log \tilde{\Omega} \left( t \right)}\\
    & \leqslant & C' \hspace*{\fill}   \text{\tmtextup{(5)}},
    \label{eq:triple.tilde.Omega}
  \end{eqnarray*}
  for some constant $C'$. Now integrating (\ref{eq:triple.tilde.Omega}) over
  $t$, we obtain a triple exponential growth estimate for $\tilde{\Omega}
  \left( t \right)$.
  
  To obtain the estimate for $\Omega \left( t \right)$, notice that
  $\tilde{\Omega} \left( t_k \right) = R \Omega \left( t_k \right)$ and
  $\tmop{loglog} x$ is a concave function for $x > e^{- 1}$, we get
  \begin{eqnarray*}
    \tmop{logloglog} \tilde{\Omega} \left( t_k \right) & = & \tmop{loglog}
    \left( \log R + \log \Omega \left( t_k \right) \right)\\
    & \leqslant & \tmop{logloglog} \Omega \left( t_k \right) + \left(
    \tmop{loglog} \right)' \left( \log \Omega \left( t_k \right) \right) \log
    R\\
    & = & \tmop{logloglog} \Omega \left( t_k \right) + \frac{\log R}{\log
    \Omega \left( t_k \right) \tmop{loglog} \Omega \left( t_k \right)} .
    \hspace*{\fill}   \text{\tmtextup{(6)}}  \label{eq:triple.log.Omega.tk}
  \end{eqnarray*}
  Combining (\ref{eq:triple.log.Omega.tk}) with the triple exponential
  estimate for $\tilde{\Omega} \left( t \right)$ and using $\Omega \left( t
  \right) < \tilde{\Omega} \left( t \right)$ for $t \in \left[ t_k, t_{k + 1}
  \right]$, we obtain (\ref{eq:triple.est.local}) immediately by taking $t =
  t_{k + 1}$.
  
  \item {\tmstrong{Global estimate.}}
  
  In the last step we obtain
  \begin{eqnarray*}
    \tmop{logloglog} \Omega \left( t_{k + 1} \right) & \leqslant &
    \tmop{logloglog} \Omega \left( t_k \right) + C \frac{R^2 r \left( 1 + C_0
    \right)}{c_L}  \left( t_{k + 1} - t_k \right)\\
    &  & + \frac{\log R}{\log \Omega \left( t_k \right) \tmop{loglog} \Omega
    \left( t_k \right)} .
  \end{eqnarray*}
  Summing over $k$ from $0$ to $n - 1$, we obtain
  \begin{eqnarray*}
    \tmop{logloglog} \Omega \left( t_n \right) & \leqslant & \tmop{logloglog}
    \Omega \left( t_0 \right) + C \frac{R^2 r \left( 1 + C_0 \right)}{c_L} 
    \left( t_n - t_0 \right)\\
    &  & \sum_{k = 0}^{n - 1} \frac{\log R}{\log \Omega \left( t_k \right)
    \tmop{loglog} \Omega \left( t_k \right)} . \hspace*{\fill}  
    \text{\tmtextup{(7)}} \label{eq:triple.log.tn}
  \end{eqnarray*}
  Now we estimate the sum in the RHS of (\ref{eq:triple.log.tn}) and prove
  that
  \begin{equation}
    \label{eq:sum.est} \sum_{k = 0}^{n - 1} \frac{\log R}{\log \Omega \left(
    t_k \right) \tmop{loglog} \Omega \left( t_k \right)} \leqslant \frac{\log
    R}{\log r} \tmop{logloglog} \Omega \left( t_n \right) + C
  \end{equation}
  for some constant $C > 0$, $n \geqslant 2$.
  
  Since $\Omega \left( t_n \right) = r^n \Omega \left( t_0 \right)$, we have
  \begin{eqnarray*}
    \sum_{k = 0}^{n - 1} \frac{\log R}{\log \Omega \left( t_k \right)
    \tmop{loglog} \Omega \left( t_k \right)} & = & \sum_{k = 0}^{n - 1}
    \frac{\log R}{\log \left( r^k \Omega_0 \right) \tmop{loglog} \left( r^k
    \Omega_0 \right)}\\
    & = & \sum_{k = 0}^{n - 1} \frac{\log R}{\left( k \log r + \log \Omega_0
    \right) \log \left( k \log r + \log \Omega_0 \right)}\\
    & = & \frac{\log R}{\log r} \sum_{k = 0}^{n - 1} \frac{\log r}{\left( k
    \log r + \log \Omega_0 \right) \log \left( k \log r + \log \Omega_0
    \right)} .
  \end{eqnarray*}
  where $\Omega_0$ is a shorthand for $\Omega \left( t_0 \right)$.
  
  Note that the sum
  \[ \sum_{k = 0}^{n - 1} \frac{\log r}{\left( k \log r + \log \Omega_0
     \right) \log \left( k \log r + \log \Omega_0 \right)} \]
  is in the form of a Riemann sum of the function $\left( x \log x \right)^{-
  1}$. This function is decreasing for $x > e^{- 1}$. Therefore the above sum
  can be bounded by
  \begin{eqnarray*}
    \int_{\log \Omega_0}^{n \log r + \log \Omega_0} \frac{1}{x \log x} \mathd
    x & = & \tmop{logloglog} \left( r^n \Omega_0 \right) - \tmop{logloglog}
    \Omega_0\\
    & = & \tmop{logloglog} \Omega \left( t_n \right) - \tmop{logloglog}
    \Omega_0
  \end{eqnarray*}
  since $\Omega \left( t_n \right) = r^n \Omega_0$ by our construction. This
  proves (\ref{eq:sum.est}).
  
  Now using the fact that $r > R$, we get
  \[ \tmop{logloglog} \Omega \left( t_n \right) \leqslant \frac{\log r}{\log r
     - \log R}  \left[ C \frac{R^2 r \left( 1 + C_0 \right)}{c_L}  \left( t_n
     - t_0 \right) + C' \right], \]
  which implies the triple exponential bound for $\Omega \left( t_n \right)$,
  and consequently no blow-up can occur at time $T_{\ast}$. This completes the
  proof of Theorem \ref{thm:triple}.
\end{itemize}

\subsubsection{Proof of Theorem \ref{thm:double}}

The proof is almost the same as that of Theorem \ref{thm:triple}. Therefore we
will only mention what are different in the four steps of the proof.
\begin{enumeratenumeric}
  \item {\tmstrong{Partition of the time interval.}} There is no difference,
  we still divide $\left[ t_0, T_{\ast} \right)$ into sub-intervals such that
  \[ \frac{\Omega \left( t_{k + 1} \right)}{\Omega \left( t_k \right)} = r \]
  for some $r > R$.
  
  \item {\tmstrong{Estimate of the lower bound of $l \left( t_k \right)$.}}
  (\ref{eq:Omega.t.kp1.t.k}) is replaced by
  \begin{equation}
    \label{eq:Omega.t.kp1.tk.double} \Omega \left( t_{k + 1} \right) \leqslant
    R \Omega \left( t_k \right)  \left[ 1 + C \frac{\left( 1 + C_0 \right)
    Rr}{c_L}  \int_{t_k}^{t_{k + 1}} \left[ \log \Omega \left( \tau \right) +
    1 \right] \mathd \tau \right] .
  \end{equation}
  \item {\tmstrong{Local triple exponential estimate.}} Define
  $\tilde{\Omega}$ in a similar way. We obtain
  \[ \left( \tmop{loglog} \tilde{\Omega} \left( t \right) \right)' \leqslant C
  \]
  for some constant $C$. Further noticing that $\tilde{\Omega} \left( t_k
  \right) = R \Omega \left( t_k \right)$, and $\log x$ is concave for $x > 0$,
  we have
  \begin{eqnarray*}
    \tmop{loglog} \tilde{\Omega} \left( t_k \right) & = & \log \left( \log R +
    \log \Omega \left( t_k \right) \right.\\
    & \leqslant & \tmop{loglog} \Omega \left( t_k \right) + \frac{\log
    R}{\log \Omega \left( t_k \right)} . \hspace*{\fill}  
    \text{\tmtextup{(8)}}
  \end{eqnarray*}
  \item {\tmstrong{Global estimate.}} From the last step we obtain
  \begin{equation}
    \tmop{loglog} \Omega \left( t_{k + 1} \right) \leqslant \tmop{loglog}
    \Omega \left( t_k \right) + C \left( t_{k + 1} - t_k \right) + \frac{\log
    R}{\log \Omega \left( t_k \right)} .
  \end{equation}
  Therefore
  \begin{equation}
    \tmop{loglog} \Omega \left( t_n \right) \leqslant \tmop{loglog} \Omega
    \left( t_0 \right) + C \left( t_{k + 1} - t_k \right) + \sum_{k = 0}^{n -
    1} \frac{\log R}{\log \Omega \left( t_k \right)} .
  \end{equation}
  Recall that $\Omega \left( t_k \right) = r^k \Omega \left( t_0 \right)$ by our
  choice of $t_k$. We have
  \begin{eqnarray*}
    \sum_{k = 0}^{n - 1} \frac{\log R}{\log \Omega \left( t_k \right)} & = &
    \frac{\log R}{\log r} \sum_{k = 0}^{n - 1} \frac{\log r}{k \log r + \log
    \Omega_0}\\
    & \leqslant & \frac{\log R}{\log r} \int_{\log \Omega_0}^{n \log r + \log
    \Omega_0} \frac{1}{x} \mathd x\\
    & = & \frac{\log R}{\log r} \left[ \tmop{loglog} \Omega \left( t_n
    \right) - \tmop{loglog} \Omega \left( t_0 \right) \right], \hspace*{\fill}
    \text{\tmtextup{(9)}}
  \end{eqnarray*}
  where the inequality is obtained by a similar argument as in the proof of
  Theorem \ref{thm:triple} using the fact that $\frac{1}{x}$ is decreasing for
  $x > 0$. Since $r > R$, we immediately obtain the double exponential
  estimate
  \begin{equation}
    \tmop{loglog} \Omega \left( t_n \right) \leqslant \frac{\log r}{\log r -
    \log R} \left[ C \left( t_n - t_0 \right) + C' \right] .
  \end{equation}
  Thus ends the proof.
\end{enumeratenumeric}

\vspace{0.2in}
\noindent
{\bf Acknowledgments.}
This work was in part supported by NSF under the NSF
FRG grant DMS-0353838 and ITR Grant ACI-0204932.

\end{document}